\documentclass[11pt,reqno]{article}
\usepackage[english]{babel}
\usepackage{amsmath,amsthm,commath,mathrsfs,amssymb,extarrows}
\usepackage{dsfont}
\usepackage{relsize}
\allowdisplaybreaks

\def \rd {\textrm{d}}

\usepackage{amsfonts}
\usepackage{graphicx}
\usepackage{enumerate}
\usepackage[right,pagewise,displaymath, mathlines]{lineno}
\usepackage{epstopdf}
\usepackage{color}
\usepackage{multirow}
\usepackage{authblk}
\usepackage[round]{natbib}
\usepackage{float}
\restylefloat{figure,table}

\usepackage{bbm}

\usepackage{authblk}
\usepackage[round]{natbib}
\usepackage{float}

\usepackage{enumitem}

\usepackage{tikz}
\usepackage{schemabloc}

\usepackage[bottom]{footmisc}

\usepackage[colorlinks = true,urlcolor = blue,citecolor = blue,breaklinks]{hyperref}

\usepackage{comment}

\usepackage{geometry}
\numberwithin{equation}{section}
\numberwithin{figure}{section}
\numberwithin{table}{section}

\newtheorem{theorem}{Theorem}[section]

\newtheorem{lemma}{Lemma}[section]

\theoremstyle{definition}

\newtheorem{remark}{Remark}[section]

\numberwithin{equation}{section}

\usepackage[title]{appendix}

\definecolor{darkred}{rgb}{0.7, 0, 0}
\definecolor{darkbrown}{rgb}{0.55, 0.2, 0.15}
\definecolor{darkblue}{rgb}{0.1,0.1,0.6}
\definecolor{darkgreen}{rgb}{0.1,0.5,0.2}

\usepackage{setspace}
\usepackage[bottom]{footmisc}

\usepackage{pgfplots}

\title{\vspace{-20mm} 
Asymptotic Expansion and Bounds for the Bias of Empirical Tail Value-at-Risk}
\author[,1,2]{Nadezhda Gribkova \thanks{e-mail: \href{mailto:n.gribkova@spbu.ru}{n.gribkova@spbu.ru}}}
\author[,3,4]{Jianxi Su \thanks{Corresponding author; e-mail: \href{mailto:jianxi@purdue.edu}{jianxi@purdue.edu}}}
\author[,4,5]{Mengqi Wang \thanks{e-mail: \href{mailto:mwan259@uwo.ca}{mwan259@uwo.ca}}}
\affil[1]{\normalsize  Saint Petersburg State University, Saint Petersburg, 199034 Russia}
\affil[2]{\normalsize  Emperor Alexander I Saint Petersburg State Transport University, \break Saint Petersburg, 190031 Russia}
\affil[3]{\normalsize Purdue University, West Lafayette, Indiana 47907, United States}
\affil[4]{\normalsize RISC Foundation, Toronto, Ontario M2N 7E9, Canada}
\affil[5]{\normalsize Western University, London, Ontario N6A 5B7, Canada}

\date{\vspace{-1mm}\normalsize\today}

\begin{document}

\maketitle

\vspace{-5mm}

{\small
\begin{center}
\textbf{Abstract}
\end{center}

Tail Value-at-Risk (TVaR) is a widely adopted risk measure playing a critically important role in both academic research and industry practice in insurance.  In data applications, TVaR is often estimated using the empirical method, owing to its simplicity and nonparametric nature.  The empirical TVaR has been explicitly advocated by regulatory authorities as a standard approach for computing TVaR. However, prior literature has pointed out that the empirical TVaR estimator is negatively biased, which can lead to a systemic underestimation of risk in finite-sample applications. This paper aims to deepen the understanding of the bias of the empirical TVaR estimator in two dimensions: its magnitude as well as the key distributional and structural determinants driving the severity of the bias.  To this end, we derive a leading-term approximation
for the bias based on its asymptotic expansion. The closed-form expression associated with the leading-term approximation enables us to obtain analytical insights into the structural properties governing the bias of the empirical TVaR estimator. To account for the discrepancy between the leading-term approximation and the true bias, we further derive an explicit upper bound for the bias.
We validate the proposed bias analysis framework via simulations and demonstrate its practical relevance using real data.
\bigskip

\noindent
{\bf Keywords}:   Expected shortfall, bootstrap, finite-sample bias, heavy-tailed distributions,  trimmed means

\medskip


\medskip

\noindent
\textbf{Dedication}:
\begin{center}
{\itshape 
This research was initiated by Ri\v{c}ardas Zitikis (1963--2025) in early 2025. This article is dedicated to his memory.
}
\end{center}
}



\newpage

\section{Introduction}
\label{intro}

The measurement of risk plays a pivotal role across a wide spectrum of applications in finance and insurance, including but not limited to, determining risk-loaded prices, assessing capital adequacy,  supporting regulatory compliance, measuring risk-based performance, and guiding portfolio optimization.   Among many risk measures employed in these contexts, Tail Value-at-Risk (TVaR), which is also referred to as AVaR, CVaR, or expected shortfall, has gained particular prominence among academics and professionals, especially following its adoption as the replacement for Value-at-Risk (VaR) as the standard risk measure within recent regulatory frameworks \citep[e.g.,][]{basel2016marketrisk,basel2019marketrisk}.  In detail, consider the risk position of a generic one-period model, denoted by $X$,  where positive values represent losses and negative values represent surpluses.  Let $F$ be the cumulative distribution function (CDF) of $X$.  The left-continuous inversion of  $F$ is given by
\begin{equation*}
\label{inversion}
F^{-1}(u)=\inf\{ x \in \mathbb{R}:~ F(x)\ge u\},\qquad u\in (0,1).
\end{equation*}
For a given probability level $p\in (0,1)$,
the TVaR of $X$ at level $p$ is defined, provided the integral is finite, as
\begin{align*}
T_p=\frac{1}{1-p}\int_p^1 F^{-1}(u) \,\rd u,
\end{align*}
which can be interpreted as the average of the upper tail quantiles of $X$.

Given the practical prominence of TVaR, a large and continually growing body of literature has been generated. A wide range of its economically grounded properties have been extensively studied, including its adherence to the coherence axioms \citep{artzner1999coherent}, comonotonic additivity \citep{kusuoka2001law}, and its representation as a Choquet integral \citep{schmeidler1989subjective} and as the solution to an expected loss minimization problem \citep{rockafellar2002conditional}.  We also refer readers to \citet{wang2021axiomatic} for an axiomatic characterization of TVaR, which also includes a comprehensive review of recent theoretical developments related to TVaR.

While different from, yet closely related to, the aforementioned literature on studying the mathematical and economic properties of TVaR, the primary objective of this paper pertains to the estimation of TVaR. 
In applied settings, the choice of estimation method is often guided by practical considerations.  Practitioners often favor approaches that are simple to implement, transparent to communicate, and whose performance does not rely on a specific set of parametric assumptions.  One such widely adopted method is the empirical, or plug-in, estimator of TVaR:
\begin{align}
\label{eqn:emperical-TVaR}
 \widehat T_{n,p}=\frac{1}{1-p}\int_p^1 \widehat F_n^{-1}(u)\,\rd u,
\end{align}
where $\widehat F_n^{-1}$ denotes the inversion of the empirical CDF $\widehat F_n$, based on $n \in \mathbb{N}$ independent and identically distributed (iid) observations $X_1,\ldots,X_n$ drawn from the distribution of $X$. 
In particular, the empirical CDF is computed via 
\begin{align*}
    \widehat F_n(x)=\frac{1}{n}\sum_{i=1}^n \mathds{1}_{\{X_i\leq x\}},
\end{align*}
where $\mathds{1}_{\{\cdot\}}$ denotes the indicator function, and its inversion can be calculated as 
\begin{align*}
    \widehat F_n^{-1}(u)=\inf\{ x \in \mathbb{R}:~ \widehat F_n(x)\ge u\} =\bigg\{
    \begin{array}{ll}
        X_{nu:n}, & \text{if $nu$ is an integer;} \\[3mm]
        X_{[nu]+1:n}, & \text{otherwise,}
    \end{array}
    \qquad u\in (0,1),
\end{align*}
where $[\cdot]$ represents the greatest integer function, and $X_{1:n}\leq\cdots\leq X_{n:n}$ are the order statistics of the samples $X_1,\dots,X_n$. 
Throughout the rest of this paper, empirical TVaR refers to the estimator in \eqref{eqn:emperical-TVaR}, unless specified otherwise.  This empirical estimator has been recommended by the European Banking Authority as the standard formula for computing TVaR \citep{eba2020rts_srm}.

The large-sample properties of empirical TVaR $\widehat T_{n,p}$ have been thoroughly studied in the related literature \citep[e.g.,][]{bjpz2008, wei2023assessing}. More general results based on the empirical tail conditional allocation have been developed in \citet{gribkova2022empirical,gribkova2022inference}. In particular, it has been shown that $\widehat T_{n,p}$ in \eqref{eqn:emperical-TVaR} is consistent.  Thereby, empirical TVaR $\widehat T_{n,p}$ ultimately provides an accurate evaluation of the true TVaR as the sample size tends to infinity.

That said, in real-world applications, data are finite and sometimes not very large.
In such cases, empirical TVaR in \eqref{eqn:emperical-TVaR} has been found to exhibit a negative bias \citep{gwz2025}. 
More precisely, existing literature has established the existence and sign of this bias, but not its magnitude. 
Consequently, two possible scenarios emerge. If the bias is small even for moderately sized samples and diminishes rapidly as the sample size increases, then empirical TVaR can still be used with confidence. On the other hand, if the bias is substantial and persistent, caution is warranted, and bias reduction techniques should be employed. Otherwise, the inherent negative bias may lead to a systematic underestimation of the riskiness of a financial position as measured by TVaR.

The goal of our paper is to examine the magnitude of the negative bias in empirical TVaR $\widehat T_{n,p}$, thereby providing insights into the two potential scenarios discussed above. While numerical methods such as jackknife \citep{quenouille1956notes,tukey1958bias} and bootstrap \citep{efron1994introduction} are commonly used to estimate the bias of an estimator for a given dataset, the resulting bias estimates are inherently tied to the specific sample at hand. As such, these methods are limited in their ability to provide a structural understanding of how the distributional characteristics of data impact the extent of the bias. Moreover, they offer little insight into how the bias evolves as the sample size increases.

In this paper, we take an analytical route to analyze the bias of empirical TVaR in \eqref{eqn:emperical-TVaR}. Specifically, following the initial conjecture the same authors of the present paper made in equation (30) of \cite{gwz2025}, this paper formally establishes the asymptotic expansion of the bias under detailed conditions and with a rigorous proof.  We then use the leading term of this expansion as a practical approximation for analyzing the bias. 
It is important to emphasize that, although the bias of the empirical TVaR estimator vanishes asymptotically due to its consistency property, the proposed leading-term approximation still offers useful practical insights, even in finite-sample settings.
This approach is analogous to using the leading term of a Taylor expansion to analyze the behavior of a complicated function. Compared with numerical resampling techniques, the closed-form expression obtained from the leading-term approximation provides explicit insight into the distributional characteristics that drive the magnitude of the bias, as well as the rate at which the bias diminishes as the sample size grows. 

The aforementioned leading-term approximation, while informative, does not generally coincide exactly with the finite-sample bias. 
Thereby, we complement it with an explicit upper bound for the negative bias. Our numerical experiments indicate that, with an appropriate choice of tuning parameters, the proposed upper bound can achieve a reasonable balance between tightness and conservativeness when assessing the bias of the empirical TVaR.

Our proposed bias analysis framework carries significant practical relevance. Beyond providing analytical insights into the structural properties of the bias of empirical TVaR, the framework delivers tractable and satisfactory approximations of the bias magnitude across a range of probability levels and sample sizes, without requiring repeated simulation or resampling procedures commonly used in numerical approaches. Understanding these relationships is particularly useful for practitioners when designing a risk assessment procedure, where trade-offs between data availability constraints and feasible or desirable choices of probability levels in TVaR applications must be balanced. Moreover, as we will see in the numerical experiments, simulation-based bias evaluation methods such as the bootstrap may, due to sampling variability, occasionally produce improper bias estimates contradicting the theoretical property that the bias of empirical TVaR is always negative. In contrast, the proposed bias analysis framework is deterministic and therefore does not suffer from this drawback. 
When a risk analyst concerns the extent of bias, the bias evaluation framework developed in this paper can be readily employed to implement bias reduction.


Finally, we note that, in addition to the empirical method considered in this paper, other commonly used approaches for estimating TVaR include parametric methods \citep[e.g.,][]{chen2023class,brazauskas2009robust,lee2010modeling}, kernel-based methods \citep[e.g.,][]{bolance2008skewed,chen2008nonparametric,yu2010kernel}, and methods originated from Extreme Value Theory \citep[e.g.,][]{goegebeur2022extreme,hua2011second,necir2010estimating}. We also refer the reader to \cite{nadarajah2014estimation} for a systematic review of different methods for estimating TVaR. It is important to note taht the objective of this paper is not to advocate for the superiority of one estimation method over another for computing TVaR. Rather, given the widespread use of the empirical TVaR estimator in \eqref{eqn:emperical-TVaR}, this paper contributes to deepening the fundamental understanding of its properties, limitations, and explores how it may be improved when necessary.

The remainder of this paper is organized as follows. Before presenting the main results, Section~\ref{sec:comment} compares the bias of the empirical TVaR estimator in \eqref{eqn:emperical-TVaR} with that of an alternative empirical estimator based on the conditional-mean formulation, which clarifies why our study focuses on the empirical TVaR estimator defined as the tail average of the quantile function. Section \ref{sec:main} then presents the proposed explicit bias analysis framework. Section \ref{sec:sim} assesses the performance of the proposed framework through an extensive simulation study. In Section \ref{sec:real-data}, the practical relevance of the proposed bias analysis framework is illustrated using the Danish fire loss data. Section \ref{sec:concl} concludes the paper. To facilitate readability, all technical proofs are relegated to Appendix \ref{app:proof}. 

\section{A commentary on empirical Tail Conditional Expectation}
\label{sec:comment}
Before presenting the proposed bias analysis framework, we devote this section to commenting on the comparison of empirical TVaR $\widehat T_{n,p}$ in \eqref{eqn:emperical-TVaR}, from the perspective of bias, with another empirical estimator based on the conditional-mean formulation.

In what follows, for notational convenience, let us shorthand the $(p\times 100)\%$-th left percentile of $F$ by
\begin{align*}
    \xi_p:= F^{-1}(p), \qquad p\in(0,1).
\end{align*}
The following relationship holds:
\begin{align*}
    T_p=\frac{1}{1-p}\int_p^1 F^{-1}(u) \,\rd u=\frac{1}{1-p}\Big[\mathbf{E}\Big(X\, \mathds{1}_{\{X>\xi_p\}}\Big) + \xi_p\times  \Big(F(\xi_p)-p\Big)\Big].
\end{align*}
Moreover, when the quantile function $F^{-1}$ is continuous around $p$, it holds that $F(\xi_p)=p$, and the above expression simplifies to
\begin{align*}
    T_p=\mathbf{E}(X\,|\, X>\xi_p),
\end{align*}
where the right-hand side corresponds to the Tail Conditional Expectation (TCE) risk functional. In this case, TVaR can be estimated using the empirical TCE estimator, defined as
\begin{equation}
\label{eqn:emperical-TCE}
\widehat {\rm TCE}_{n,p}=\frac{1}{n-[np]}\sum_{i=[np]+1}^n X_{i:n}.
\end{equation}

We argue that the empirical TVaR estimator of interest, namely $\widehat T_{n,p}$ in \eqref{eqn:emperical-TVaR}, exhibits a smaller magnitude of negative bias (i.e., a smaller absolute bias) than $\widehat {\rm TCE}_{n,p}$ in \eqref{eqn:emperical-TCE}, regardless of whether the distribution $F$ is continuous or not. Therefore, from the perspective of bias minimization, $\widehat T_{n,p}$ in \eqref{eqn:emperical-TVaR} should be preferred over $\widehat {\rm TCE}_{n,p}$ in \eqref{eqn:emperical-TCE} for computing TVaR. To see the above claim, let us rewrite empirical TVaR in terms of the following more computational friendly form: 
\begin{align}
\label{connection_1}
\widehat{T}_{n,p}=\frac{1}{1-p}\int_p^1\widehat{F}_n^{-1}(u)\,\rd u&=\frac{1}{1-p}\left[ \int_p^{\frac{[np]+1}{n}} \widehat{F}_n^{-1}(u)\,\rd u+\int_{\frac{[np]+1}{n}} ^1 \widehat{F}_n^{-1}(u)\,\rd u\right]\notag\\
&=\frac{[np]+1-np}{(1-p)n} \, X_{[np]+1:n}+\frac{1}{(1-p)n}\sum_{i=[np]+2}^nX_{i:n}\notag\\
&=-\frac{np-[np]}{(1-p)n}\, X_{[np]+1:n}+\frac{1}{(1-p)n}\, \sum_{i=[np]+1}^nX_{i:n}\notag\\
&=-\frac{np-[np]}{(1-p)n}\, X_{[np]+1:n}+\frac{n-[np]}{(1-p)n}\, \widehat{\rm TCE}_{n,p}.
\end{align}
Accordingly, we have
\begin{align}
\label{connection}
\notag
\widehat{\rm TCE}_{n,p}&=\frac{(1-p)n}{n-[np]}\, \widehat T_{n,p}+\frac{np-[np]}{n-[np]}X_{[np]+1:n}\\
&=\widehat T_{n,p}-\frac{np-[np]}{n-[np]}\Big(\widehat T_{n,p}- X_{[np]+1:n}\Big).
\end{align}
Recall that the bias of $\widehat{T}_{n,p}$ is always negative \citep{gwz2025}.  Moreover, note that
\begin{align*}
    \widehat{T}_{n,p}=\frac{1}{1-p}\int_p^1\widehat{F}^{-1}_n(u)\,\rd u \geq \frac{1}{1-p}\int_p^1\widehat{F}^{-1}_n(p)\,\rd u = \widehat{F}^{-1}_n(p).
\end{align*}
If $np$ is not an integer, then $\widehat{F}_n^{-1}(p)=X_{[np]+1:n}$, hence $\widehat{T}_{n,p}\geq X_{[np]+1:n}$.
Thereby, the subtraction term in \eqref{connection} is always non-negative, which introduces an additional negative bias to the estimator $\widehat{\rm TCE}_{n,p}$ compared to $\widehat T_{n,p}$, unless $np$ is an integer, in which case this additional bias is zero.  In other words, the empirical estimator of interest $\widehat T_{n,p}$ always outperforms $\widehat{\rm TCE}_{n,p}$ in the sense of having a smaller, or at the very least an equal, negative bias. This observation justifies our focus on $\widehat T_{n,p}$ as the estimator of interest in this paper.

\section{Main results
}
\label{sec:main}

We aim to study the bias of empirical TVaR $\widehat T_{n,p}$ in \eqref{eqn:emperical-TVaR} under minimal assumptions, so that the results can be broadly applicable across a wide range of settings. Nevertheless, certain finiteness conditions are still required to guarantee that the problem under study is well-posed. The first condition ensures that the risk measure of interest, $T_p$, is well-defined: For a given $p\in (0,1)$
\begin{enumerate}[label={\rm(C$_{\arabic*}$)}]
\item \label{c1}
\qquad\qquad\qquad\qquad\qquad \qquad \qquad$\displaystyle \int_p^1 F^{-1}(u)\, \rd u<\infty$.
\end{enumerate}
In addition, we need the mean of the empirical TVaR estimator, and thus its bias, to exist. This is ensured by the following condition: For a given $p\in (0,1)$,
\begin{enumerate}[label={\rm(C$_{\arabic*}$)}]
\setcounter{enumi}{1}
\item \label{c2}
\qquad\qquad\qquad\qquad\qquad $\displaystyle \mathbf{E}\Big(|X|^{\varepsilon}\, \mathds{1}_{\{X\leq \xi_p\}}\Big)<\infty$, for some $\displaystyle \varepsilon>0$.
\end{enumerate}
It is worth noting that condition \ref{c2} is very mild. In particular, it may fail only in the presence of a sufficiently heavy left tail. In the context of risk management, loss models are typically assumed to have positive support or at the very least lower bounded. In these cases, condition \ref{c2} can be automatically satisfied.


 To understand why conditions \ref{c1} and \ref{c2} together imply the finiteness of the bias,
we shall examine the difference between $\widehat{T}_{n,p}$ and $T_p$.  To this end, we need to introduce some additional notations.
For $i=1,\ldots,n$, let $W_i$ denote the winsorized counterpart of $X_i$ with respect to the threshold $\xi_p$:
\begin{equation*}
\label{W_i}
W_i=\max(X_i,\, \xi_p)=\begin{cases}\xi_p,&X_i\leq \xi_p;\\
X_i,& X_i>\xi_p.
\end{cases}
\end{equation*}
Note that $W_1,\ldots,W_n$ are iid random variables. We denote their common mean by 
\begin{align*}
    \mu_{\scriptscriptstyle W}=\mathbf{E}(W_1).
\end{align*}
Moreover, let $N_p:={\rm {\rm card}}\{i\,:\, X_i \leq \xi_p\}$ denote the number of iid samples $X_1,\ldots,X_n$ that fall below or at the threshold $\xi_p$. 

The following lemma spells out $(\widehat{T}_{n,p}-T_p)$ explicitly in terms of a sum of iid centered winsorized random variables, accompanied by a non-positive remainder term. At the outset, let us define
\begin{equation}
\label{notations_1}
a\wedge b:=\min(a,b);\quad a\vee b:=\max(a,b);\quad {\rm {\rm {\rm {\rm sign}}}}(a):=\begin{cases}\frac {a}{|a|},&a\neq 0,\\
0,&a=0.\end{cases}
\end{equation}

\begin{lemma}
\label{lamma_1} Suppose that condition \ref{c1} holds, then we have
\begin{equation}
\label{main_rep}
\widehat{T}_{n,p}-T_p=\frac 1{n(1-p)}\sum_{i=1}^n\big(W_i-\mu_{\scriptscriptstyle W}\big)-R_n,
\end{equation}
where
\begin{align}
\label{eqn:residual}
R_n=\frac{np-[np]}{n(1-p)}(X_{[np]+1:n}-\xi_p) - \frac {{\rm {\rm {\rm {\rm sign}}}}(N_p-[np])}{n(1-p)}\sum_{i=([np]\wedge N_p)+1}^{[np]\vee N_p} (X_{i:n}-\xi_p).
\end{align}
Furthermore,
\begin{align}
\label{R_n_positive}
R_n\geq 0.
\end{align}
\end{lemma}

From here on, let us denote the bias of empirical TVaR $\widehat T_{n,p}$ by
\begin{align*}
    B_n=\mathbf{E}\big(\widehat{T}_{n,p}\big)-T_p.
\end{align*} 
From the representation \eqref{main_rep} in Lemma \ref{lamma_1}, it follows that 
\begin{align}
        B_n=-\mathbf{E}(R_n).
        \label{eqn:B-R}
\end{align}
Therefore, the finiteness of the bias $B_n$ is equivalent to the finiteness of $\mathbf{E}(R_n)$, or equivalently, to the finiteness of the expectations of the order statistics involved in the expression of $R_n$ in \eqref{eqn:residual}.  By Proposition 2 of \cite{s1974}, conditions \ref{c1} and \ref{c2} together ensure the existence of moments of central order statistics of any order, and in particular, the finiteness of the first moment, which is sufficient for our purposes in the current context.
\begin{remark}
\label{rmk:why-c}
\itshape
By the relationships noted in \eqref{R_n_positive} in Lemma \ref{lamma_1} as well as in \eqref{eqn:B-R}, it follows that the bias $B_n$ is non-positive, which is in agreement with the result established in \cite{gwz2025}.
\end{remark}



We are now in a position to present one of the main results of this paper, which can offer analytical insights into the structural properties of the bias $B_n$ through the leading term of its asymptotic expansion. 

\begin{theorem}\label{Thm_1}
For a fixed $p\in (0,1)$, suppose that conditions \ref{c1} and \ref{c2} hold.
Further, assume that the probability density function (PDF) of 
$X$, denoted by $f$, exists and is continuous in a neighborhood of the point $\xi_p=F^{-1}(p)$, and that $f(\xi_p)>0$.
Then, the following relationship holds for the bias of empirical TVaR $\widehat T_{n,p}$: 
\begin{align}
\label{B_n}
B_n=\ell_n+o\big(n^{-1}\big),\qquad \ell_n=-\frac{p}{2nf(\xi_p)}.
\end{align}
\end{theorem}


\begin{remark}
\itshape
{In equation (30) of \cite{gwz2025}, the same authors as in the present paper conjectured the asymptotic expansion in \eqref{B_n}, motivated by heuristic reasoning derived from Edgeworth-type expansions for the distribution of normalized trimmed means studied in \cite{gh2006,gh2007}. That conjecture was stated without a complete specification of the required conditions nor a rigorous proof. Theorem~\ref{Thm_1} of the present paper formally confirms the correctness of this initial conjecture by establishing the expansion under a complete and explicit set of conditions.
}
\end{remark}

Thanks to the explicit expression established in Theorem \ref{Thm_1}, several analytical insights emerge immediately, which are summarized in the succeeding remarks.
\begin{remark}
\itshape
    The asymptotic expansion in \eqref{B_n} shows that the negative bias decreases at the rate of $n^{-1}$ as the sample size increases.  It also confirms that the empirical TVaR estimator is asymptotically unbiased.
\end{remark}

\begin{remark}
\label{rmk:varying_p}
\itshape
For a fixed sample size $n$, the leading term $\ell_n$ in \eqref{B_n} is proportional to $p/f(\xi_p)$. Assuming that $f$ decreases monotonically along the far right tail, which is typically observed in empirical (insurance) data, then the function $p\mapsto p/f(\xi_p)$ is increasing for sufficiently large values of $p$. This implies that the negative bias becomes more pronounced as the probability level $p$ increases. Intuitively, increasing $p$ reduces the size of the tail region and leads to fewer observations above the VaR threshold $\xi_p$, so the estimation of the tail mean becomes more difficult and its finite-sample bias becomes larger.
\end{remark} 

\begin{remark}
\label{rmk:varying_a}
\itshape
When $n$ and $p$ are fixed, the magnitude of the leading term $\ell_n$ is inversely related to $f(\xi_p)$. In risk management practice, TVaR is commonly evaluated at high probability levels. 
When a data distribution becomes more heavy-tailed, its associated probability density tends to spread farther into the extreme region. As a result, the density around a fixed high quantile becomes smaller, or equivalently $f(\xi_p)$ tends to decrease.
While this opposite relationship between tail heaviness and $f(\xi_p)$ may not hold universally for all distribution families, it is indeed the case for models commonly used in risk modeling practice.  For instance, the Pareto distribution considered in our succeeding simulation study in Section \ref{sec:sim} statisfies this relationship. Consequently, heavier tails lead to a larger value of the term $1/f(\xi_p)$, and thus a more pronounced negative bias of the empirical TVaR $\widehat T_{n,p}$. 
\end{remark}

In addition, we note that the relationship between the leading-term bias of empirical TVaR and $f(\xi_p)$ stands in stark contrast with kernel-based TVaR estimators, whose leading bias term is increasing with $f(\xi_p)$ rather than decreasing \citep[e.g.,][]{chen2008nonparametric,handkernel}. The reasoning behind this observation may be that  
a smaller density at $\xi_p$ implies fewer observations and larger spacing between upper order statistics in the neighborhood surrounding $\xi_p$. As a result, an empirical estimator is less capable of accurately locating the true quantile threshold, which leads to a more pronounced downward bias. However, for kernel methods, a kernel estimator smooths the loss distribution in a neighborhood of the quantile level and therefore assigns positive weight to observations just below $\xi_p$. These sub-threshold observations are smaller than the true tail losses, and their inclusion pulls the TVaR estimate downward. When $f(\xi_p)$ is larger, a greater amount of probability mass accumulates near $\xi_p$, amplifying this smoothing-induced downward pull and can increase the magnitude of the negative bias.


We note that in practical applications,  the quantile $\xi_p$ and the PDF $f$ are unknown. To evaluate the leading term $\ell_n$ in \eqref{B_n}, one may estimate $\xi_p$ using the empirical quantile and estimate $f(\xi_p)$ using a kernel density estimator or another nonparametric method, such as local likelihood density estimation, local polynomial density estimation, or an orthogonal series estimator, and so forth. 

Thus far, Theorem \ref{Thm_1} identifies the leading term in the asymptotic expansion of the bias $B_n$ of the empirical TVaR. Although informative, the leading term $\ell_n$ need not match the finite-sample bias exactly. 
To strengthen our analysis of the bias $B_n$, we next complement the leading-term approximation in Theorem \ref{Thm_1} with an explicit upper bound for $|B_n|$, which yields a more conservative evaluation of the bias. 


\begin{theorem}\label{Thm_2}
For a fixed $p\in (0,1)$, suppose that conditions \ref{c1} and \ref{c2} are satisfied. Assume also that
the inversion $F^{-1}$ is locally H\"{o}lder continuous of order $\gamma\in (0,1]$ at the
point $p$, that is, there exists a neighborhood $\mathbb{U}_p$ of the point $p$ and a constant
$c_{\gamma}>0$ such that
\begin{equation}\label{Lipschitz}
|F^{-1}(u)-F^{-1}(p)|
\le c_{\gamma}\,|u-p|^{\gamma},
\quad u\in\mathbb{U}_p.
\end{equation}
Then, for any $\delta>0$, there exists $n_0$ such that the bias satisfies
\begin{equation}\label{B_n_estim}
-B_n
\;\le\;
\frac{C\,p^{\frac{1+\gamma}{2}}}
{n^{\frac{1+\gamma}{2}}(1-p)^{\frac{1-\gamma}{2}}},
\qquad n\ge n_0,
\end{equation}
where $C=c_\gamma\, (1+\delta)$.
\end{theorem}

\begin{remark}
\label{rmk:holder}
\itshape
The local H\"{o}lder continuity constant $c_\gamma$ in \eqref{Lipschitz} quantifies the local smoothness of the quantile function $F^{-1}$ near 
$p$. A larger value of $c_\gamma$ indicates greater steepness or variability of 
$F^{-1}$ in a neighborhood of $p$.  Hence, a larger value of $c_\gamma$ implies that even small fluctuations in the sample 
order statistics may translate into notable errors in the empirical TVaR estimate. 
Consequently, a larger upper bound is needed to accommodate the increased sensitivity 
of the negative bias, which explains why the bias bound in \eqref{B_n_estim} is monotonically increasing in 
$c_\gamma$ when the other parameters are held constant.  

Moreover, since the derivative of the upper bound with respect to $c_\gamma$ is 
proportional to $n^{-(1+\gamma)/2}$, its sensitivity to $c_\gamma$ decreases as the 
sample size $n$ grows. This implies that for small samples, an accurate assessment of 
$c_\gamma$ has a much stronger impact on correctly implementing the upper bound than it does for large 
samples. We will later discuss how to estimate $c_\gamma$ in practice. 
\end{remark}

\begin{remark}
\label{rmk:slack}
\itshape
In the choice of $C = c_\gamma(1+\delta)$, the slack parameter $\delta>0$ is not a structural parameter of the data distribution but 
an arbitrarily small inflation factor introduced to ensure that the upper bound 
\eqref{B_n_estim} holds uniformly for all sufficiently large $n$. A smaller $\delta$ yields a tighter upper bound in \eqref{B_n_estim}, but the sample size $n$ may then be larger before the inequality is guaranteed. 
\end{remark}

As discussed in Remark \ref{rmk:holder}, implementing the negative-bias bound in 
Theorem \ref{Thm_2} requires an appropriate evaluation of the local H\"{o}lder continuity 
constant $c_\gamma$. First, let us consider an ideal scenario where the underlying data distribution $F$ is 
known.  For almost all commonly used continuous parametric 
distributions with a positive and continuous density at $\xi_p$, the H\"{o}lder exponent of $F^{-1}$ near $p$ satisfies 
$\gamma=1$. To calculate $c_1$, namely the local H\"{o}lder constant $c_{\gamma}$ when $\gamma=1$, note that for a small interval $\mathbb{X}_{\xi_p,h} = [\xi_p-h,\xi_p+h]$, $h>0$, the mean value 
theorem implies that for all $u$ sufficiently close to $p$,
\[
|F^{-1}(u)-F^{-1}(p)| \le 
\Big(\sup_{x\in \mathbb{X}_{\xi_p,h}} \frac{1}{f(x)}\Big)\,|u-p|.
\]
Thereby, a local H\"{o}lder  (or Lipschitz in this case of $\gamma=1$) constant can be computed via 
\begin{align}
c_1 = \Bigl(\inf_{x\in \mathbb{X}_{\xi_p,h}} f(x)\Bigr)^{-1}.
\label{eqn:c_1_cal}
\end{align}

In a more realistic scenario where the data distribution is unknown and only data are available, a natural estimator of $c_1$ is obtained by replacing $f$ with a 
local density estimator near $\widehat\xi_p$, the empirical estimate of the $(p\times 100)\%$-th quantile. Let $\widehat{f}(x)$ denote such an estimator 
(e.g., kernel density estimator) and let 
$\mathbb{X}_{\widehat\xi_p,h}=[\widehat\xi_p-h,\, \widehat\xi_p+h]$ be a small data-driven neighborhood. 
Then an empirical Lipschitz constant can be
\begin{align}
\widehat{c}_1 = 
\Big(\inf_{x\in \mathbb{X}_{\widehat \xi_p,h}} \widehat{f}(x)\Big)^{-1}.
\label{eqn:c_1_est}
\end{align}

\section{Simulation analysis}
\label{sec:sim}
The simulation analysis in this section serves two complementary purposes. The first is to assess how well the explicit leading bias identified in Theorem~\ref{Thm_1} reflects the behavior of the finite-sample bias across varying sample sizes, probability levels, and tail behaviors. The second is to evaluate the effectiveness of the leading-term approximation from Theorem~\ref{Thm_1} and the finite-sample bound from Theorem~\ref{Thm_2} in approximating the true bias of empirical TVaR $\widehat T_{n,p}$.

Throughout this simulation study, we assume that the risk random variable $X$ follows a Pareto distribution with CDF:
\begin{align}
    F(x) = 1 - x^{-\alpha}, \qquad x > 1,
    \label{eqn:Pareto-cdf}
\end{align}
where $\alpha > 0$ denotes the shape parameter. The smaller the value of $\alpha$, the heavier the tail of the distribution of $X$. The mean of $X$ is finite when $\alpha>1$, and so is the TVaR $T_p$ for any $p\in(0,1)$.  This is the case of interest throughout this section. 

\subsection{Leading-term approximation}
Let us begin by examining the leading-term approximation of the bias of empirical TVaR.
Under the model given in \eqref{eqn:Pareto-cdf}, the mean of the $i$-th order statistic $X_{i:n}$ based on an iid sample $X_1,\ldots,X_n$ from $X$ can be computed explicitly as, for $i < n+1-\alpha^{-1}$,
\begin{align*}
    \mathbf{E}\big(X_{i:n}\big)
    = \frac{\Gamma(n+1)\,\Gamma\bigl(n-i+1 - {\alpha}^{-1}\bigr)}
           {\Gamma(n-i+1)\,\Gamma\bigl(n+1 - {\alpha}^{-1}\bigr)},
\end{align*}
where $\Gamma(\cdot)$ represents the gamma function.  Consequently, using the representation in the third line of \eqref{connection_1}, the exact bias of empirical TVaR can be derived explicitly as
\begin{align*}
    B_n&=\mathbf{E}\Big[-\frac{np-[np]}{(1-p)n}\, X_{[np]+1:n}+\frac{1}{(1-p)n}\, \sum_{i=[np]+1}^nX_{i:n}\Big] - T_p\\
    &=\frac{\Gamma(n+1)}{(1-p)\, n\,\Gamma\bigl(n+1-\alpha^{-1}\bigr)}
\left[
-(np-[np])\,
\frac{\Gamma\bigl([np]+1-\alpha^{-1}\bigr)}{\Gamma([np]+1)}
+
\sum_{i=[np]+1}^n
\frac{\Gamma\bigl(i-\alpha^{-1}\bigr)}{\Gamma(i)}
\right]
- T_p,
\end{align*}
where the true TVaR of \eqref{eqn:Pareto-cdf} is given by
\begin{align*}
    T_p=\frac{\alpha}{\alpha-1}(1-p)^{-1/\alpha},\qquad p\in(0,1).
\end{align*}

We consider the following baseline parameter setting: $n = 500$, $p = 95\%$, and $\alpha = 3$. We then perturb each parameter individually to examine how the leading-term bias $\ell_n$ varies with $n$, $p$, and $\alpha$, and to verify that the resulting behavior of the bias of empirical TVaR aligns with the analytical insights established in Section~\ref{sec:main}.
Two versions of the leading-term bias $\ell_n$ are evaluated. The first assumes that the data-generating model \eqref{eqn:Pareto-cdf} is known, in which case $\ell_n$ can be computed exactly. The second scenario considers a more realistic situation in which only the simulated samples are available.  In this case, the kernel density estimator with a Gaussian kernel is used to estimate the unknown PDF appearing in $\ell_n$.
Because of the sampling variability inherent in the simulated data, each configuration of the perturbed parameters is replicated $m = 100$ times to assess variability in the resulting empirical leading-term bias estimates.

\begin{figure}[h!]
    \centering
    \includegraphics[width=0.5\linewidth]{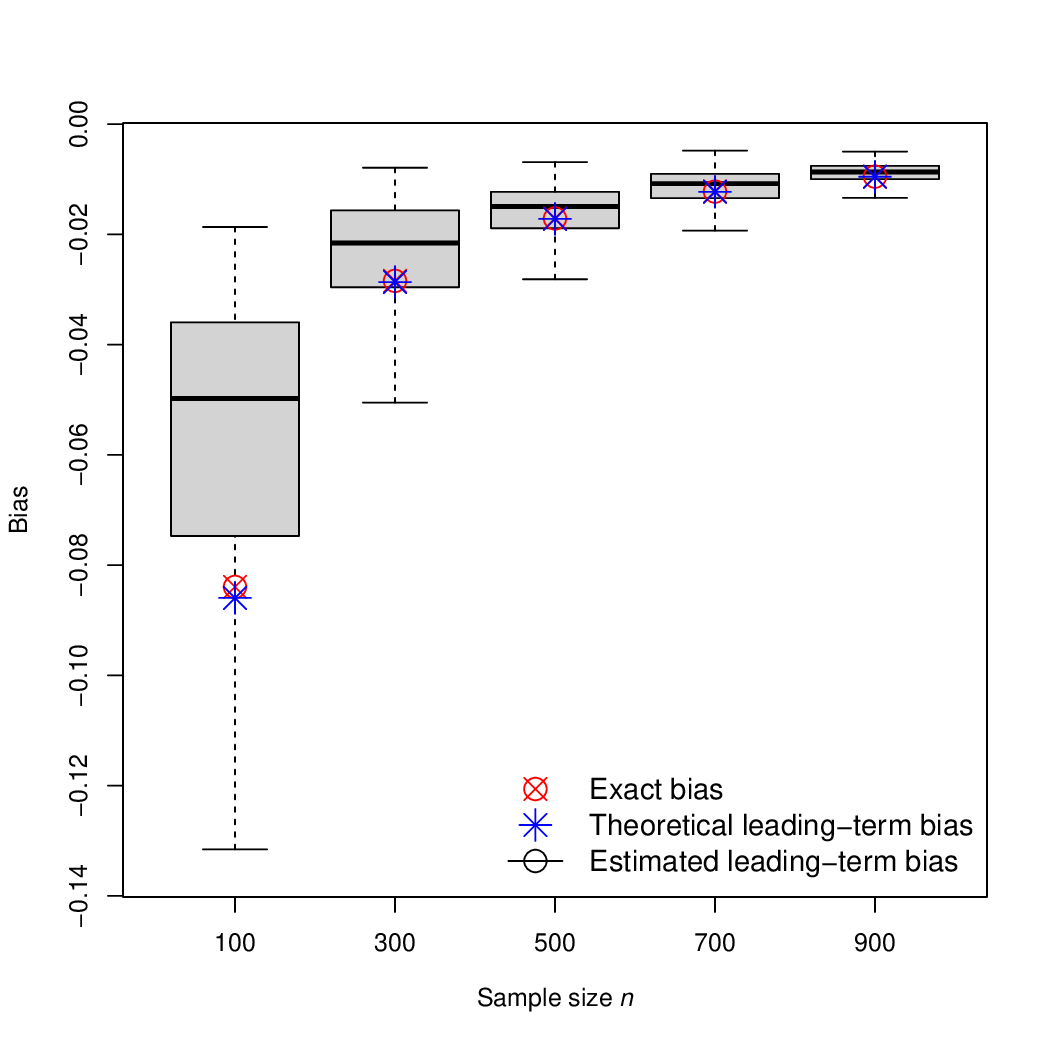}
    \caption{Comparison of the exact TVaR bias (red circle-cross markers), the theoretical leading-term approximation (blue star markers), and the estimated leading-term approximation (boxplots), with fixed parameters $p=95\%$, $\alpha=5$, and varying sample sizes $n\in \{100,\, 300,\, 500,\, 700,\, 900\}$.}
    \label{fig:bias_n}
\end{figure}

Figure \ref{fig:bias_n} presents a comparison of the exact TVaR bias with the leading-term bias $\ell_n$ derived in Theorem \ref{Thm_1} with varying sample sizes. In the figure, the boxplots are constructed from $m = 100$ simulated datasets of size $n$. As shown, when the data-generating distribution is assumed to be known, the exact bias and its leading-term approximation agree very closely. Among the sample sizes considered, a visible yet minor discrepancy is observed only when $n = 100$. This indicates that, at least for the Pareto model, the residual term in the asymptotic expansion \eqref{B_n} is negligible even at small sample sizes.
When the distribution is assumed to be unknown and $f(\xi_p)$ must be estimated from data, the performance of the leading-term approximation deteriorates due to the limited accuracy of the kernel-based estimation of $f(\xi_p)$. For smaller samples (e.g., $n=100$ or $300$), the exact bias $B_n$ can deviate substantially from the estimated leading-term bias, occasionally falling near or outside the whiskers of the boxplots. However, as the sample size increases, the estimation of $f(\xi_p)$ becomes more accurate and stable. Thereby, for larger sample sizes (e.g., $n=900$), the leading-term approximation based on the kernel density estimate tracks the exact bias of empirical TVaR very closely.

\begin{figure}[h]
    \centering
    \includegraphics[width=0.5\linewidth]{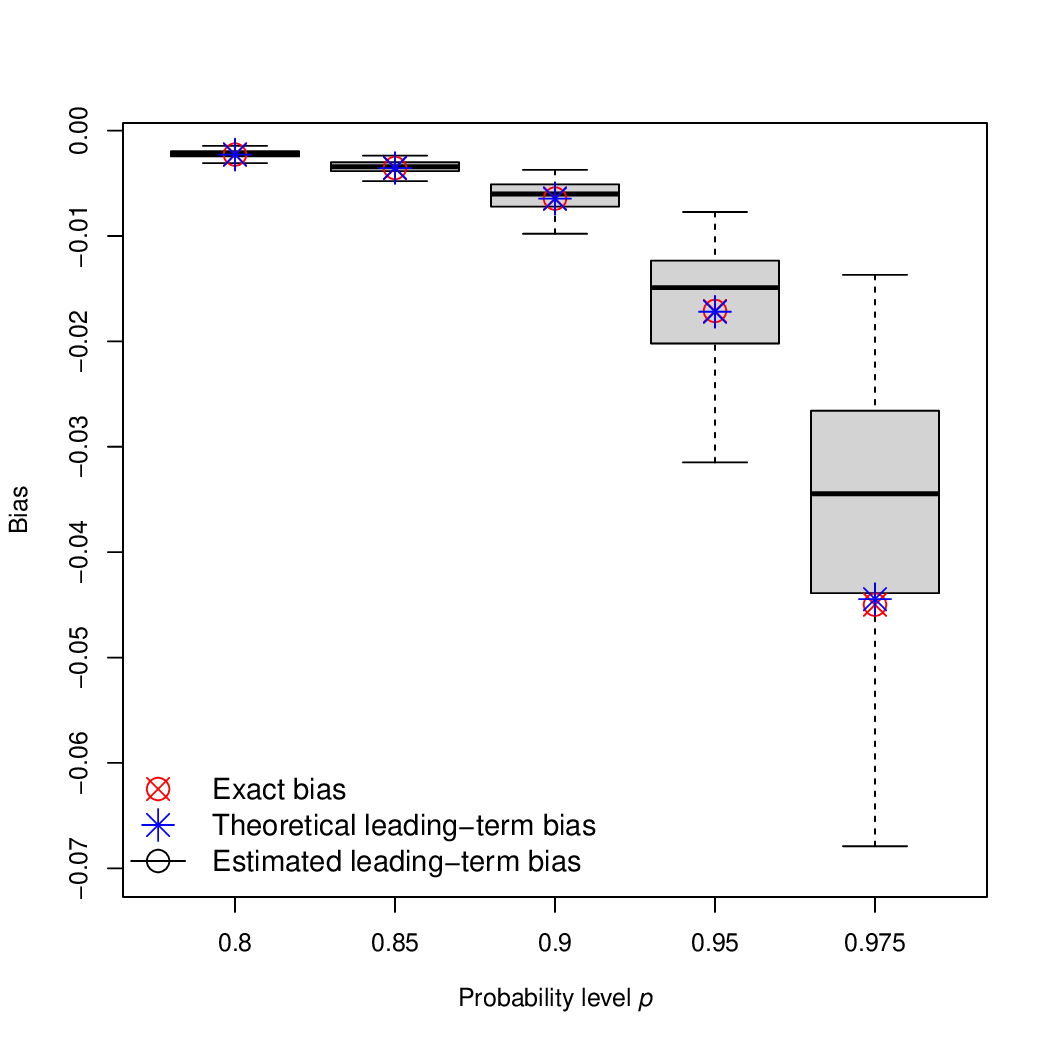}
    \caption{Comparison of the exact TVaR bias (red circle-cross markers), the theoretical leading-term approximation (blue star markers), and the estimated leading-term approximation (boxplots), with fixed parameters $n=500$, $\alpha=5$, and varying probability levels $p\in \{80,\, 85,\, 90,\, 95,\, 97.5\}\%$.}
    \label{fig:bias_p}
\end{figure}

Next, we examine the impact of the probability level $p$ on both the exact bias $B_n$ and the leading-term bias $\ell_n$. As $p$ increases, fewer observations fall beyond the quantile threshold $\xi_p$, thereby reducing the amount of data available for estimating TVaR. Consequently, the negative bias of empirical TVaR estimator becomes more pronounced. This pattern observed in Figure~\ref{fig:bias_p} is consistent with the analytical insights discussed in Remark~\ref{rmk:varying_p}. 
When the data-generating distribution is treated to be unknown, increasing $p$ also deteriorates the accuracy of the kernel density estimator for $f(\xi_p)$, leading to a larger discrepancy between the estimated leading-term bias and the exact bias. Nevertheless, the theoretical leading-term bias tracks the exact bias very closely across all probability levels considered. This suggests that the discrepancy observed in the estimated leading-term bias is mainly attributable to the estimation error in the kernel density estimator, rather than caused by the asymptotic expansion in \eqref{B_n}.

\begin{figure}[h]
    \centering
    \includegraphics[width=0.5\linewidth]{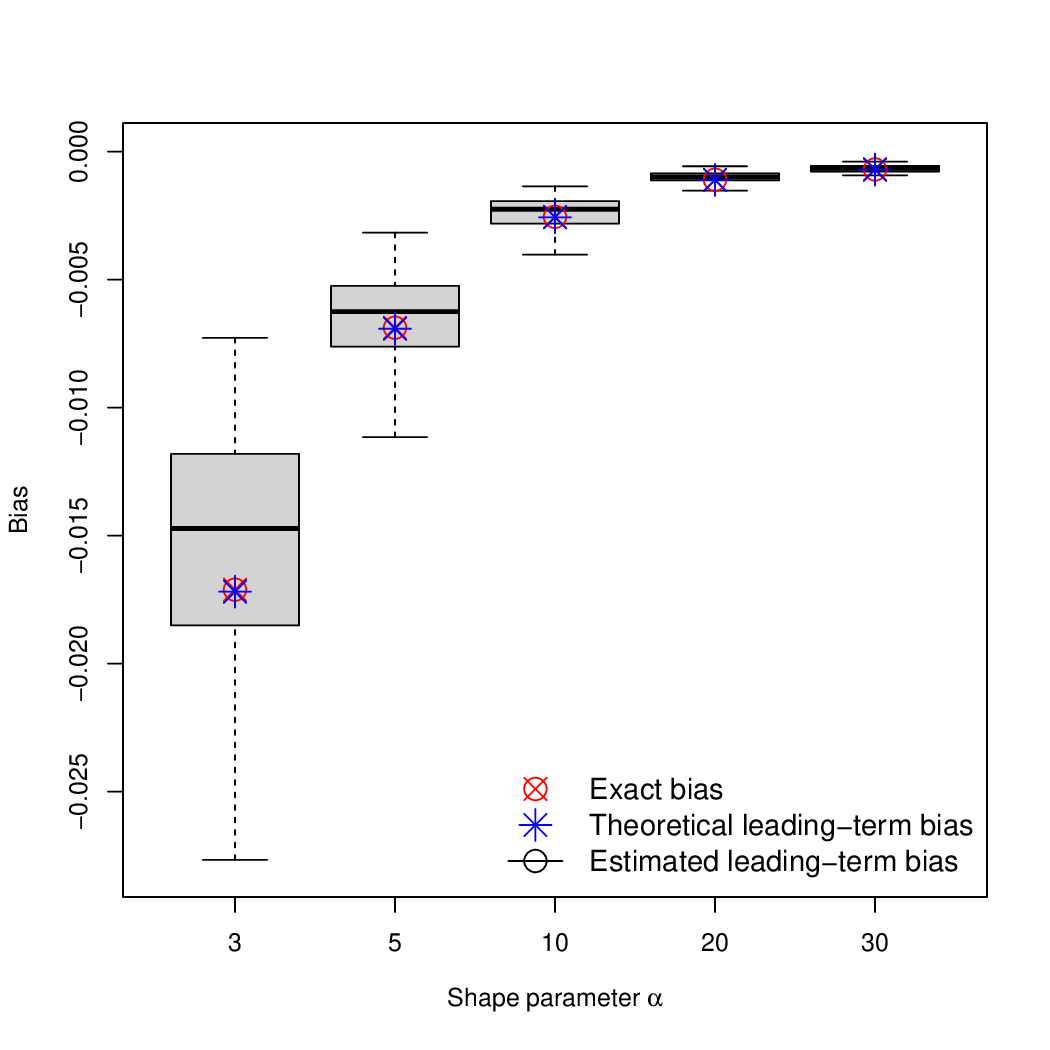}
    \caption{Comparison of the exact TVaR bias (red circle-cross markers), the theoretical leading-term approximation (blue star markers), and the estimated leading-term approximation (boxplots), with fixed parameters $n=500$, $p=95\%$, and varying values of the shape parameter $\alpha\in \{3,\, 5,\, 10,\, 20,\, 30\}$.}
    \label{fig:bias_a}
\end{figure}

Further, we study the impact of the tail heaviness of the data distribution in \eqref{eqn:Pareto-cdf}, controlled by the shape parameter $\alpha$, on the bias of empirical TVaR $B_n$ and its leading-term approximation $\ell_n$. The results are summarized in Figure  \ref{fig:bias_a}. As shown, smaller values of $\alpha$, corresponding to heavier tails, lead to a larger magnitude of the negative bias for both the exact TVaR and its leading-term approximation. This observation is consistent with the analytical insights discussed in Remark \ref{rmk:varying_a}. In particular, under data model in \eqref{eqn:Pareto-cdf}, elementary algebraic calculations yield  
\begin{align*}
f(\xi_p) = \alpha\,(1-p)^{(\alpha+1)/\alpha},
\end{align*}
which is increasing in $\alpha>0$ for a fixed $p \in (0,1)$. Thus, smaller values of $\alpha$, thus heavier tails, lead to smaller values of $f(\xi_p)$, which in turn amplify the negative leading-term bias.  This behavior of $\ell_n$ matches that of the exact bias $B_n$ across all values of $\alpha$ considered. We also observe from Figure~\ref{fig:bias_a} that heavier tails impair the accuracy of estimating $f(\xi_p)$ via kernel density estimation, thus resulting in a more noticeable discrepancy between the estimated leading-term bias and the exact bias when $\alpha$ is small.

\subsection{Negative-bias upper bound}
We now turn to examining how well the negative-bias upper bound in Theorem \ref{Thm_2} captures the true bias of empirical TVaR. To implement the bound, it is necessary to specify two tuning parameters. The first one is the slack parameter $\delta$, which controls the trade-off between the tightness of the bound and the sample size required for its validity in finite samples (see Remark~\ref{rmk:slack} for a detailed discussion). The second one is the locality parameter $h$, which determines the local Lipschitz constant through \eqref{eqn:c_1_cal}. Although introduced from different perspectives, note that these two tuning parameters in fact play similar roles in practice. Specifically, a larger value of the locality parameter $h$ enlarges the interval over which the Lipschitz constant is evaluated, resulting in a larger estimate of $c_1$ and, consequently, a more conservative upper bound. This effect parallels the role of the slack parameter $\delta$ discussed in Remark \ref{rmk:slack}. For this reason, throughout the analysis in this subsection, we fix the locality parameter $h = 0.05$ and focus on examining the sensitivity of the upper bound with respect to different choices of the slack parameter $\delta$.

\begin{figure}[h!]
    \centering
\includegraphics[width=0.5\linewidth]{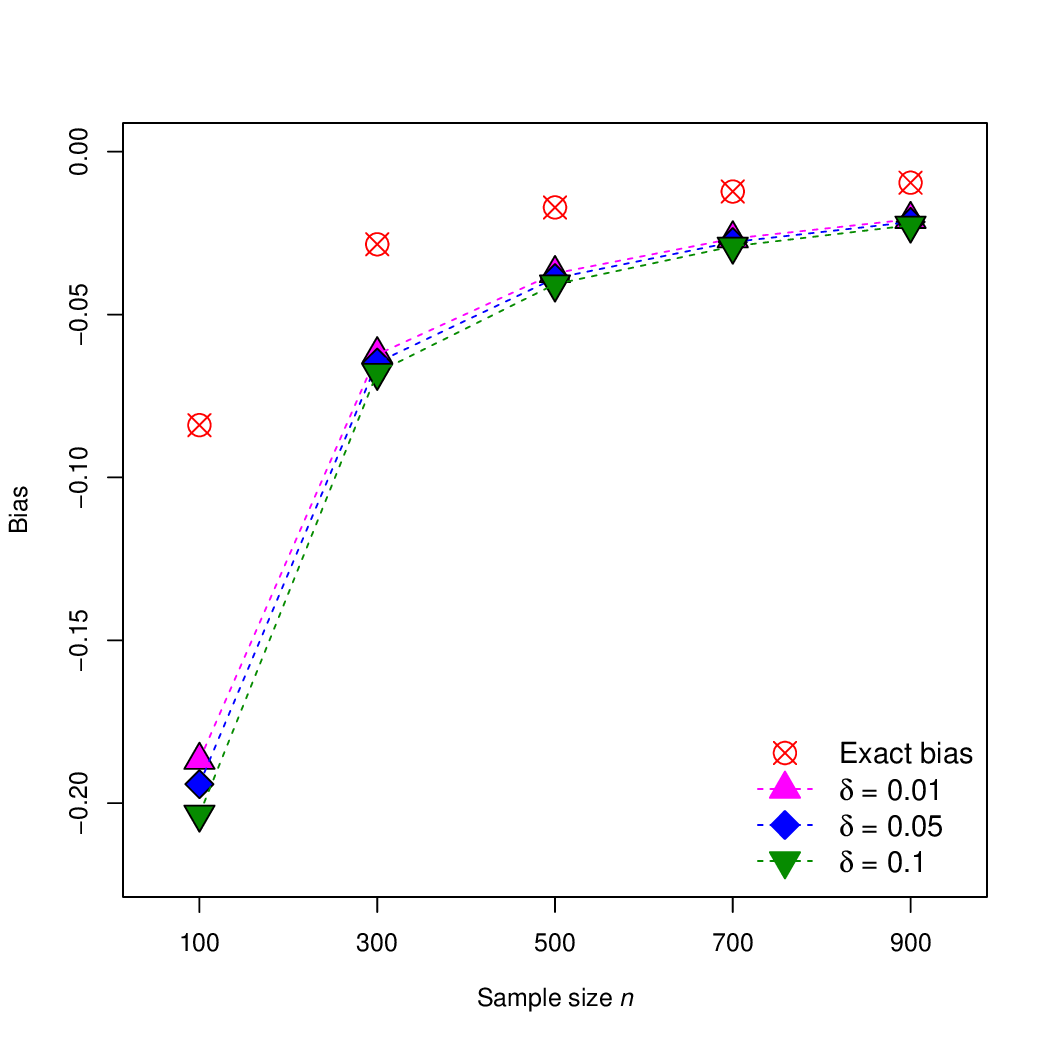}
\caption{Comparison of the exact TVaR bias with the theoretical negative-bias upper bound, with fixed parameters $p = 95\%$, $\alpha = 3$, and $h = 0.05$, across varying sample sizes $n \in \{100,\,300,\,500,\,700,\,900\}$ and different choices of the slack parameter $\delta \in \{0.01,\,0.05,\,0.1\}$.
}
    \label{fig:bound_n}
\end{figure}

Suppose that the data-generating model \eqref{eqn:Pareto-cdf} is known. For a given locality parameter $h$, the local Lipschitz constant defined in \eqref{eqn:c_1_cal} admits the following closed-form expression:
\begin{align*}
c_1
&=
\left(\inf_{x\in\mathbb{X}_{\xi_p,h}} f(x)\right)^{-1}
=
\frac{(\xi_p+h)^{\alpha+1}}{\alpha}
=
\frac{\bigl((1-p)^{-1/\alpha}+h\bigr)^{\alpha+1}}{\alpha},
\end{align*}
where the second equality follows from the fact that the PDF of the Pareto distribution in \eqref{eqn:Pareto-cdf} is strictly decreasing.
We vary the sample size $n$ to examine whether the negative-bias upper bound remains valid for small sample sizes and, if so, how tight the bound is across different values of $n$. The comparison results are summarized in Figure \ref{fig:bound_n}.

As can be observed in Figure \ref{fig:bound_n}, the derived upper bound always covers the true negative bias across all values of the slack parameter $\delta$ and sample sizes $n$ considered. For a fixed $\delta$, the bound captures the decreasing trend of the exact negative bias as the sample size $n$ increases. For a fixed $n$, larger values of the slack parameter $\delta$ lead to a more conservative bound and thus a larger discrepancy from the exact bias, which is consistent with the discussion in Remark \ref{rmk:slack}. Interestingly, the gap between the exact bias and the bound is substantially larger for smaller sample sizes and lessens noticeably as $n$ increases. This behavior can be attributed to the presence of positive residual terms in the derived upper bound under the Pareto model \eqref{eqn:Pareto-cdf}, which reduce as the sample size grows. As a result, the bound becomes increasingly tight in larger samples.

\begin{figure}[h!]
    \centering
\includegraphics[width=0.5\linewidth]{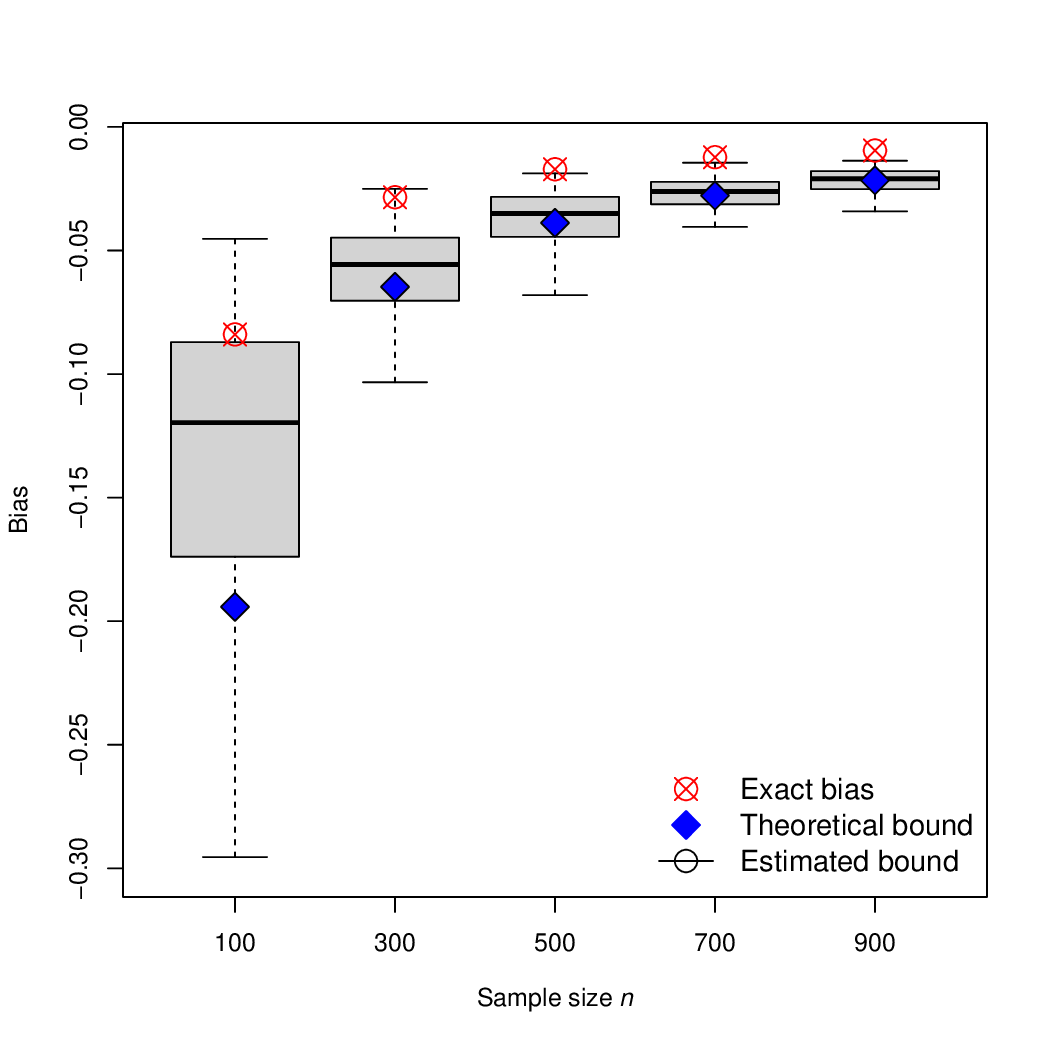}
\caption{Comparison of the exact TVaR bias (red circle-cross markers), the theoretical negative-bias upper bound (blue diamond markers), and the estimated negative-bias upper bound (boxplots), with fixed parameters $p = 95\%$, $\alpha = 3$, $\delta = 0.05$, and $h = 0.05$, across varying sample sizes $n \in \{100,\,300,\,500,\,700,\,900\}$.
}
    \label{fig:est_bound_n}
\end{figure}

In another scenario where the data-generating model \eqref{eqn:Pareto-cdf} is treated as unknown, the kernel-based estimator \eqref{eqn:c_1_est} is employed to estimate the Lipschitz constant required for computing the negative-bias upper bound. Similar to numerical experiments conducted in the previous subsection, each configuration is repeated $m=100$ times to account for the sampling variability inherent in the simulated data.  
The resulting performance of the bound relative to the exact negative bias is summarized in Figure \ref{fig:est_bound_n}. Several observations can be made. First, the boxplots corresponding to the estimated upper bound always cover the theoretical upper bound. When the sample size is small, the kernel density estimator more frequently misestimates the PDF in the tail region, leading to a noticeable discrepancy between the center of the boxplots and the theoretical upper bound. As a result, the estimated upper bound may occasionally fail to capture the exact bias. As the sample size increases, the kernel density estimator becomes more accurate and stable. As a result, the boxplots become narrower, and their centers converges toward the theoretical upper bound.  In this example, we observe that once the sample size exceeds 500, the estimated upper bound consistently covers the true bias, which indicates improved finite-sample reliability of the negative-bias upper bound as the sample size grows.

\section{Real-data illustration}
\label{sec:real-data}
We now consider a real-data illustration of the theoretical results obtained in this current paper.  Our real data example is based on the Danish fire loss dataset \citep{mcneil1997estimating}. The data consist of 2\,167 individual fire insurance losses recorded between 1980 and 1990, which have been inflation-adjusted to 1985 price levels and are reported in millions of Danish kroner. 
Prior studies of the dataset have documented that its distribution exhibits a heavy right tail, with a tail index of approximately 1.4 \citep[e.g.,][]{mcneil1997estimating,resnick1997discussion}. As demonstrated by the simulation experiments in Section \ref{sec:sim}, heavy-tailed patterns can induce an increased bias in $\widehat{T}_{n,p}$, particularly at high probability levels. This motivates us to examine the magnitude of the bias of the estimator.

Although the tail of the distribution of the Danish fire data is relatively heavy, as indicated by the estimated tail index of approximately 1.4, the distribution still has a finite first moment. Consequently, Condition \ref{c1} is satisfied. Moreover, since the loss data are strictly positive, Condition \ref{c2} also holds. These conditions together justify the validity of using the leading-term bias approximation in Theorem \ref{Thm_1}, as well as the corresponding upper bound in Theorem \ref{Thm_2}, to analyze the bias of empirical TVaR.

In the absence of knowledge of the true underlying loss distribution, unlike in the simulation study, the exact bias is not available for direct comparison. To address this, we compare the proposed leading-term approximation and upper bound against estimates obtained using a classical numerical approach, such as the bootstrap method. Specifically, we employ the nonparametric Efron bootstrap with 
$
l=1\,000$ resample sets, each of size 
$
n=2\,167$, which is the same as the size of the original dataset. To account for the variability inherent in the bootstrap procedure, we repeat the bootstrap bias estimation 
$
m=100$ times.

\begin{figure}[h!]
    \centering
\includegraphics[width=0.5\linewidth]{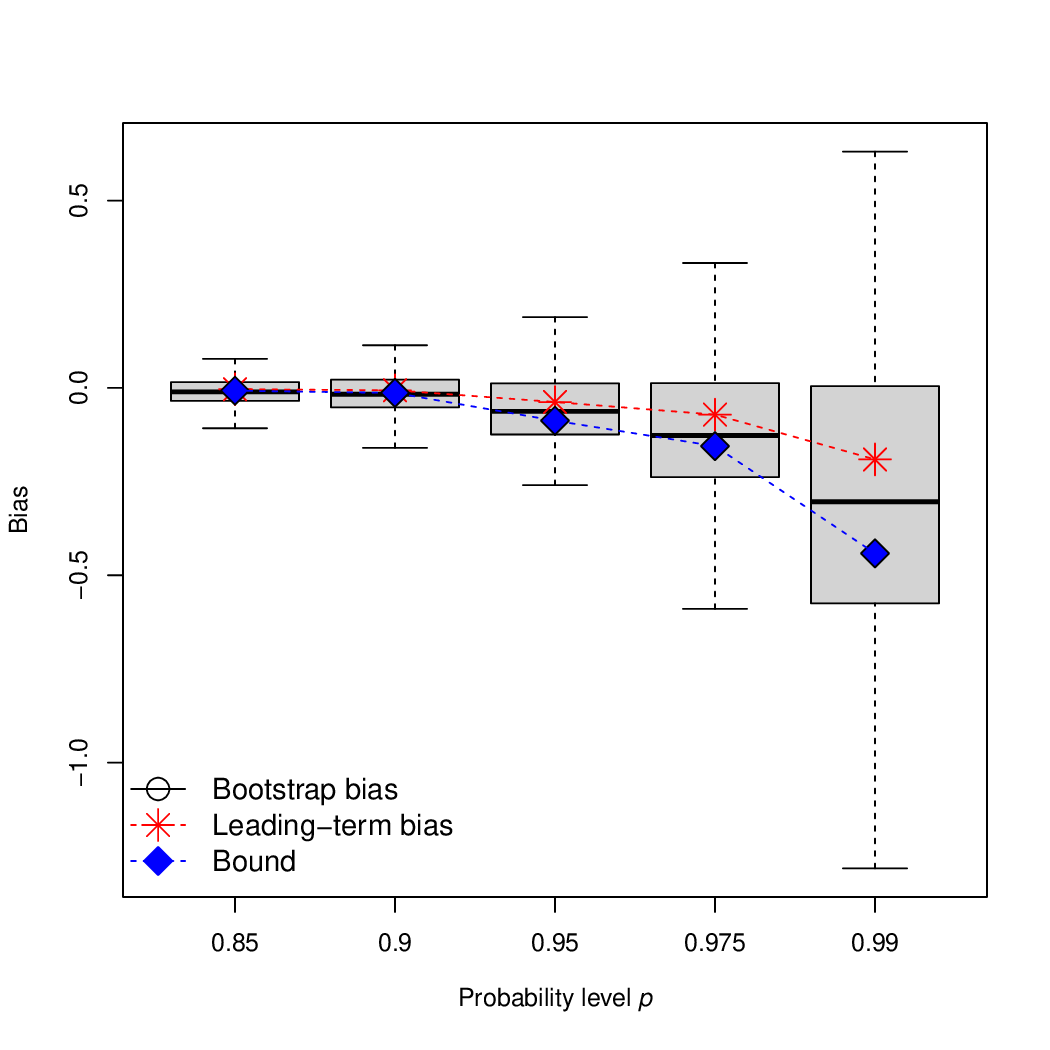}
\caption{A comparison of bias estimates for the empirical TVaR based on the Danish fire loss dataset, obtained via the leading-term approximation (red star markers), the negative-bias upper bound (blue diamond markers), and the bootstrap method (boxplots).
}
    \label{fig:fire_bias_p}
\end{figure}

Figure \ref{fig:fire_bias_p} depicts a comparison of bias estimates for the empirical TVaR \eqref{eqn:emperical-TVaR} obtained using different methods, including the leading-term approximation and the negative-bias upper bound derived in Theorems \ref{Thm_1} and \ref{Thm_2}, respectively, as well as the bootstrap method. We observe that, due to resampling variability, the bootstrap procedure may occasionally produce positive bias estimates. Such estimates contradict the theoretical result that the bias of empirical TVaR $B_n$ is always negative \citep{gwz2025}; also see discussion in Remark \ref{rmk:why-c}. For the two proposed approximations, owing to their analytical nature, they do not suffer from this drawback induced by resampling variation, and therefore provide bias assessments that are always consistent with the negativity property.

At lower probability levels, the leading-term approximation matches closely with the means of the bootstrap estimates. At higher probability levels, the leading-term approximation exhibits more pronounced discrepancies from the mean of the bootstrap estimates. At the same time, the bootstrap estimates themselves display increased variability, which is caused by the greater sampling uncertainty associated with observations in the extreme tail region that drive the estimation of TVaR. Moreover, across all probability levels considered, the proposed negative-bias upper bound consistently covers both the leading-term approximation and the means of the bootstrap bias estimates, while remaining reasonably close to the lower edge of the bootstrap boxplot whiskers. This behavior indicates a satisfactory performance of the negative-bias upper bound in terms of achieving a practical balance between tightness and conservativeness.

Compared with the bootstrap approach which is a numerical method, the proposed analytical methods including the leading-term bias approximation and the negative-bias upper bound, are more convenient to implement, as they do not rely on repeated resampling. In addition to the analytical insights already discussed in Section \ref{sec:main}, these methods can provide practitioners with an efficient and transparent framework for assessing the significance of bias. In particular, they support informed decision-making regarding whether bias-reduction strategies are warranted or whether additional data collection is desirable. 

For instance, in this data example, $\widehat{T}_{n,p}$ yields a TVaR estimate of approximately $24$ at the $95\%$ probability level, while the corresponding bias of the empirical TVaR is estimated to be around $-0.1$, representing less than $0.05\%$ of the TVaR value. This is perhaps an affirmative finding, as such a negligible level of bias indicates that empirical TVaR can be used with confidence even without applying any bias correction. For more sophisticated users seeking to further reduce the bias, the leading-term approximation $\ell_n$ can be readily applied for bias reduction. Moreover, when TVaR is used for economic capital calculation, a conservative analyst may adopt the negative-bias upper bound to incorporate an additional buffer that accounts for the bias inherent in the TVaR estimation process.

\begin{figure}[h!]
    \centering
\includegraphics[width=0.45\linewidth]{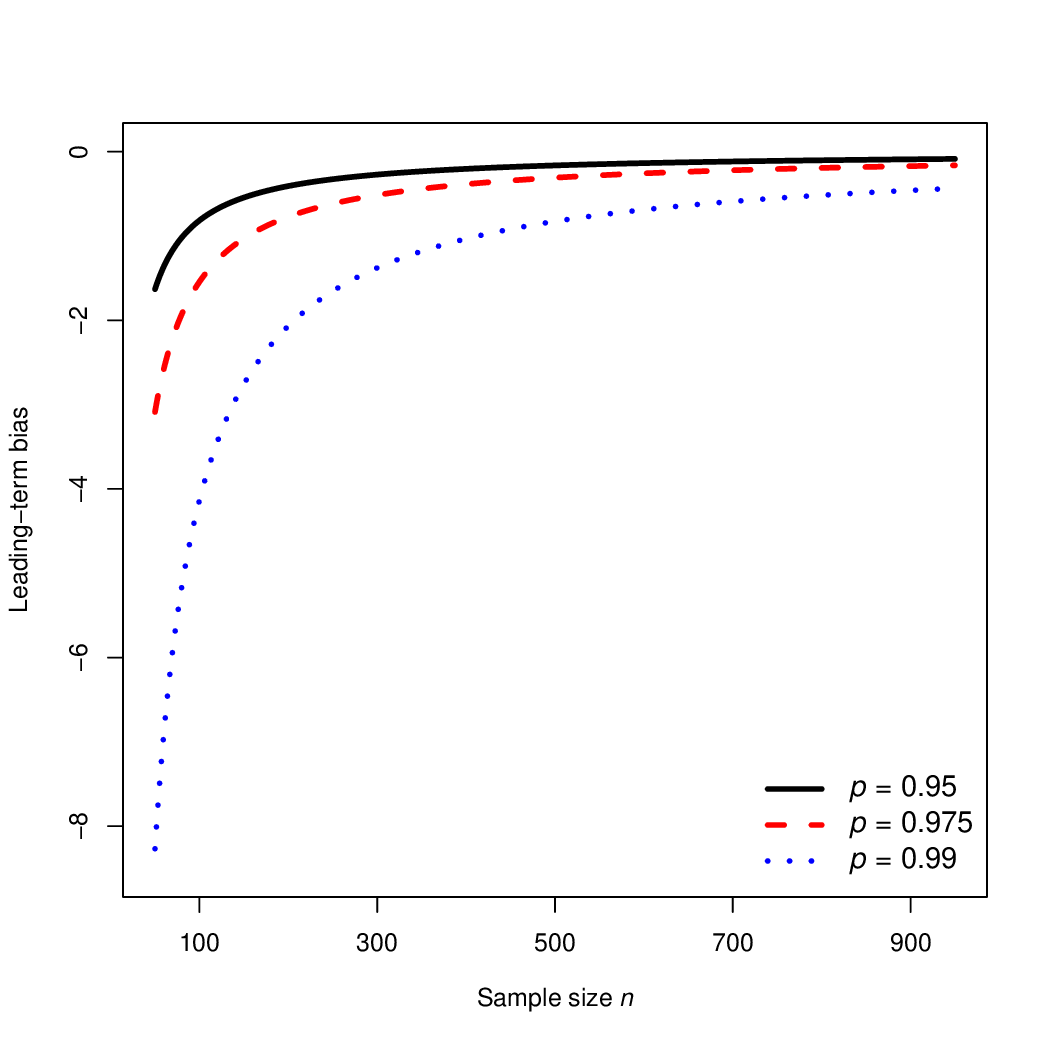}
\includegraphics[width=0.45\linewidth]{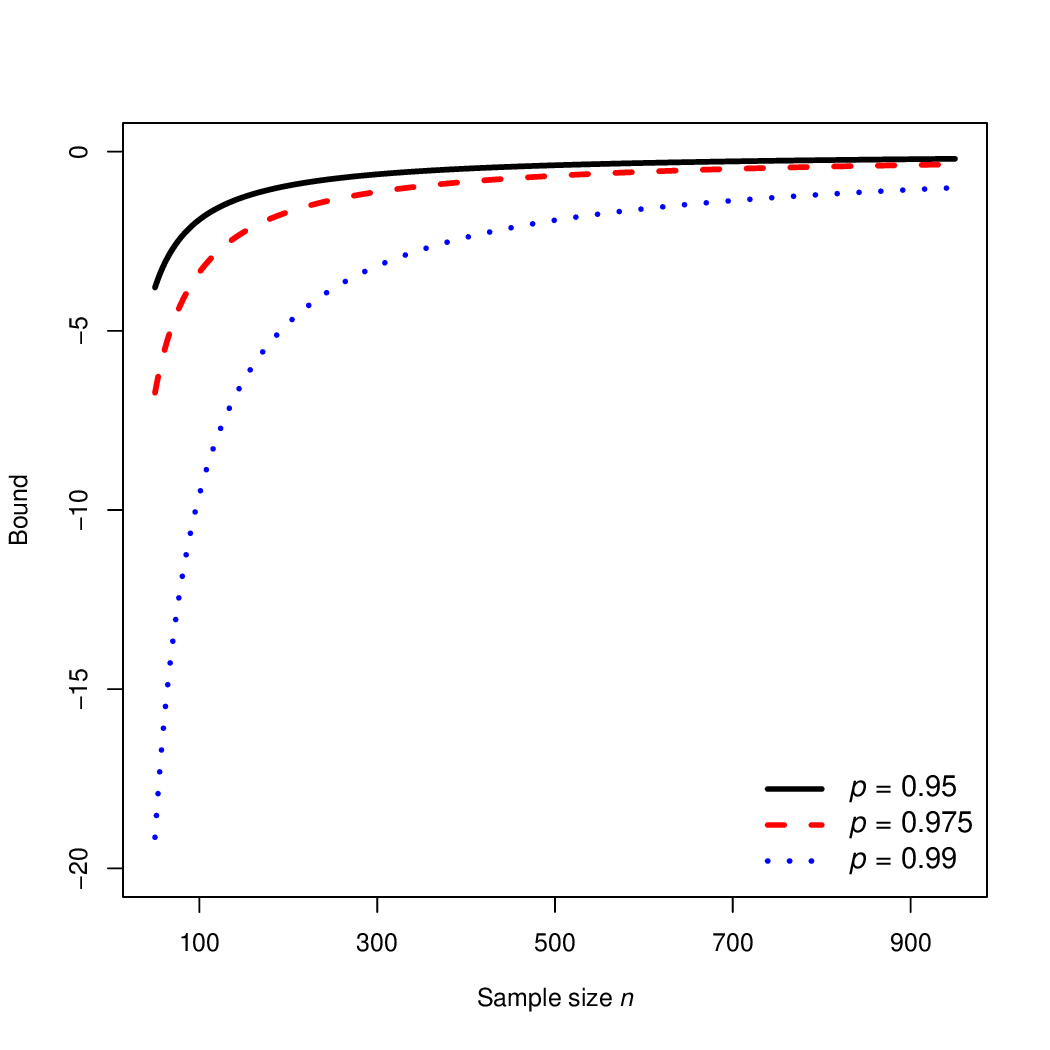}
\caption{Leading-term bias (left panel) and negative-bias upper bound (right panel) of the empirical TVaR as functions of the sample size $n$, evaluated at probability levels $p\in \{95,\, 97.5,\, 99 \}\%.$
}
    \label{fig:fire_bias_curves}
\end{figure}

For analysts interested in understanding how different choices of sample size and probability level affect the accuracy of TVaR estimation, the derived leading-term approximation and the negative-bias upper bound provide useful tools for assessing these relationships. In Figure \ref{fig:fire_bias_curves}, we first use the original Danish fire loss data to estimate $f(\xi_p)$ and the Lipschitz constant $c_1$, which are required for implementing the leading-term approximation $\ell_n$ in \eqref{B_n} and the negative-bias upper bound in \eqref{B_n_estim}. We then plot the resulting bias versus sample-size curves across different probability levels. These curves can offer insights into whether the marginal benefit of acquiring additional data is meaningful, or whether adjusting the probability level used in the calculation of TVaR is more appropriate given the current sample size.

\section{Conclusions}
\label{sec:concl}
TVaR has been widely adopted by practitioners and regulators as a conservative risk measure for capturing tail risk. In practical applications, TVaR is typically estimated using the empirical method due to its nonparametric nature and ease of implementation. It is known that the empirical TVaR estimator exhibits a negative finite-sample bias \citep{gwz2025}. However, despite its practical importance, the magnitude of this bias and its dependence on the distributional characteristics of the data and the sample size remain underexplored.

This paper addresses this gap by developing two analytical approaches for evaluating the bias of the empirical TVaR. The first approach is based on the leading term in the asymptotic expansion of the bias, while the second approach establishes an explicit upper bound on the negative bias. These results yield explicit analytical insights into how data distributional properties, sample sizes, and probability levels jointly determine the magnitude of the bias. Simulation studies demonstrate that the proposed methods perform satisfactorily in approximating the exact bias, even in moderately small sample settings. We further apply the proposed approaches to study the bias involved in TVaR estimation using the Danish fire loss data, and demonstrate how the resulting insights can inform prospective study design and risk assessment procedures, especially when balancing data availability constraints against feasible and desirable choices of probability levels in TVaR applications.   

\section*{Acknowledgment}
Mengqi Wang's research has been supported by the NSERC Alliance-MITACS Accelerate grant (ALLRP 580632-22) entitled ``New Order of Risk Management: Theory and Applications in the Era of Systemic Risk'' from the Natural Sciences and Engineering Research Council (NSERC) of Canada and the national research organization Mathematics of Information Technology and Complex Systems (MITACS) of Canada.

\bibliography{ref}
\bibliographystyle{apalike}

\begin{appendices}
\section{Technical proofs}
\label{app:proof}
Let us introduce some additional notation to facilitate the succeeding technical developments. Recall that the winsorized version of $X_i$ with respect to the threshold $\xi_p$
is denoted by  $W_i=\max(X_i,\, \xi_p)$, $i=1,\ldots,n$. The random variables $W_1,\ldots,W_n$ are iid.
The common CDF of $W_1,\ldots,W_n$ can be computed via
\begin{equation*}
\label{F_w}
F_{\scriptscriptstyle W}(x)=
\begin{cases}
0, \ & x < \xi_p,\\
F(x), \ & x\geq \xi_p,
\end{cases}
\end{equation*}
and the corresponding quantile function is given by 
\begin{equation*}
\label{F_W^{-1}}
F_{\scriptscriptstyle W}^{-1}(u)=\max(\xi_p,F^{-1}(u))=
\begin{cases}
\xi_p, \ &  u\leq p,\\
F^{-1}(u), \ & u >p.
\end{cases}
\end{equation*}
The common mean of $W_1,\ldots,W_n$, shorthanded by $\mu_{\scriptscriptstyle W}=\mathbf{E}(W_1)$, can be expressed as
\begin{equation}
\label{E_W}
\mu_{\scriptscriptstyle W}=\int_0^1F^{-1}_{\scriptscriptstyle W}(u)\,\rd u=p\,\xi_p+\int_p^1F^{-1}(u)\, \rd u.
\end{equation}

Moreover, let $W_{1:n}\leq \cdots \leq W_{n:n}$ denote the order statistics based on the random variables $W_1,\dots,W_n$. For $i=1,\ldots,n$, it holds that
\begin{equation*}
\label{W_{i:n}}
W_{i:n}=
\begin{cases}
\xi_p, \ & i \leq N_p,\\
X_{i:n}, \ & i > N_p.
\end{cases}
\end{equation*}

\begin{proof}[\bf Proof of Lemma~\ref{lamma_1}] Recall the representation given in the third line of \eqref{connection_1}.  Multiplying both sides by $(1-p)$, we get:
\begin{equation*}
\label{rep_1}
(1-p)\widehat{T}_{n,p}=-\frac{np-[np]}{n}X_{[np]+1:n}+\frac{1}{n}\sum_{i=[np]+1}^nX_{i:n}.
\end{equation*}
Consider the average of $W_1,\ldots,W_n$.  We have
\begin{equation*}
\label{rep_2}
\frac 1n\sum_{i=1}^nW_i=\frac 1n\sum_{i=1}^nW_{i:n}=\frac{N_p}{n}\xi_p+\frac 1n\sum_{i=N_p+1}^n X_{i:n},
\end{equation*}
because $W_{i:n}=X_{i:n}$ for all $i>N_p$. Taking into account \eqref{E_W}, we find
\begin{align}
\label{rep_3}
\notag
(1-p)\big(\widehat{T}_{n,p}-T_p\big)-\frac 1n\sum_{i=1}^n(W_i-\mu_{\scriptscriptstyle W})&=-\frac{np-[np]}{n}X_{[np]+1:n}+\frac{1}{n}\sum_{i=[np]+1}^nX_{i:n}-\frac{N_p}{n}\xi_p-\frac 1n\sum_{i=N_p+1}^nX_{i:n}\\[2mm]
\notag
&\quad -\int_p^1F^{-1}(u)\,\rd u+\underbrace{p\,\xi_p+\int_p^1F^{-1}(u)\,\rd u}_{\mu_{\scriptscriptstyle W}}\\[2mm]
&=-\frac{np-[np]}{n}\big(X_{[np]+1:n}-\xi_p\big)-\frac{N_p-[np]}{n}\xi_p+S_n,
\end{align}
where
\begin{equation}
\label{rep_4}
S_n=\frac 1n\sum_{i=[np]+1}^nX_{i:n}-\frac 1n\sum_{i=N_p+1}^n X_{i:n}=\frac {{\rm {\rm {\rm {\rm sign}}}}(N_p-[np])}n\sum_{i=([np]\wedge N_p)+1}^{[np]\vee N_p} X_{i:n}.
\end{equation}

Combining \eqref{rep_3} and \eqref{rep_4}, we obtain 
\begin{align*}
\label{rep_5}
\notag
\widehat{T}_p-T_p-\frac {1}{n(1-p)}&\sum_{i=1}^n(W_i-\mu_{\scriptscriptstyle W})\\
= &-\frac{np-[np]}{n(1-p)}\big(X_{[np]+1:n}-\xi_p\big) + \frac {{\rm {\rm {\rm {\rm sign}}}}(N_p-[np])}{n(1-p)}\sum_{i=([np]\wedge N_p)+1}^{[np]\vee N_p} (X_{i:n}-\xi_p)\\
=&-R_n,
\end{align*}
which implies \eqref{main_rep} and \eqref{eqn:residual}.

Next, we prove \eqref{R_n_positive}. There are three cases: (i) $[np]=N_p$; (ii) $[np]<N_p$; and (iii)~$[np]>N_p$. For case (i), we have $X_{[np]+1:n} - \xi_p>0$, hence the first term of $R_n$ is not negative, whereas its second term is zero because, according to the definition of the sign function in \eqref{notations_1}, ${\rm {\rm {\rm {\rm sign}}}}(N_p-[np])=0$. 

For case (ii),  we have $X_{i:n}\leq \xi_p$ for $i=[np]+1,\dots, N_p$, and hence,
\begin{align*}
R_n&=\frac{np-[np]}{n(1-p)}\big(X_{[np]+1:n}-\xi_p\big) - \frac {1}{n(1-p)}\sum_{i=[np]+1}^{N_p} (X_{i:n}-\xi_p)\\
&=\frac{np-[np]-1}{n(1-p)}\big(X_{[np]+1:n}-\xi_p\big)- \frac {\mathds{1}_{\{N_p>[np]+1\}}}{n(1-p)}\sum_{i=[np]+2}^{N_p} (X_{i:n}-\xi_p)\geq 0.
\end{align*}
Here, it is noteworthy that the indicator $\mathds{1}_{\{N_p>[np]+1\}}$ in front of the second sum is needed because if $N_p=[np]+1$, then the second sum disappears.

For case (iii), we have ${\rm {\rm {\rm {\rm sign}}}}(N_p-[np])=-1$, and $X_{i:n}>\xi_p$ for all $i>N_p$. Hence,
\begin{equation*}
R_n=\frac{np-[np]}{n(1-p)}\big(X_{[np]+1:n}-\xi_p\big) + \frac {1}{n(1-p)}\sum_{i=N_p+1}^{[np]} (X_{i:n}-\xi_p)\geq 0.
\end{equation*}

The proof for this assertion is now finished.
\end{proof}

\bigskip

To establish the main theorems in Section \ref{sec:main}, the following assertion serves as an important auxiliary result. 
Specifically, let $k$ and $\varepsilon$ be arbitrary positive numbers, and set $\rho := k/\varepsilon$. 
The following lemma suggests an upper bound estimate for the absolute moment of the
$i$-th order statistic, $\mathbf{E}\lvert X_{i:n}\rvert^{k}$, for
$\rho \le i \le n-\rho+1$, in terms of the absolute moment of the
underlying distribution, $\mathbf{E}\lvert X\rvert^{\varepsilon}$.
 The proof of the lemma can be found in \cite{g1995}.  

\begin{lemma}[\citeauthor{g1995}, \citeyear{g1995}]
\label{estim_mom}
For all $n \ge 2\rho + 1$ and all $i$ such that $\rho \le i \le n - \rho + 1$, the following inequality holds:
{\em
\begin{align}
\label{L_es_mom1}
\mathbf{E}\lvert X_{i:n}\rvert^{k} 
    \le 
    C(\rho)\left\{
        \mathbf{E}\lvert X\rvert^{\varepsilon}
        \, h\!\left(\frac{i}{\,n+1\,}\right)
    \right\}^{\rho},
\end{align}}
where the function $h:(0,1)\to(0,\infty)$ is given by $h(u)=[{u(1-u)]}^{-1}$, $u\in(0,1)$, 
and the constant $C(\rho)$ may be taken as
$
C(\rho)=2\sqrt{\rho}\,\exp\!\left(\rho+{7}/{6}\right).
$
\end{lemma}

\bigskip




To facilitate the presentation of the proofs of Theorems \ref{Thm_1} and \ref{Thm_2}, 
some additional notations are needed. Let $U_1,\dots,U_n$ be independent $[0,1]$-uniform random variables, and $U_{1:n}\leq\cdots\leq U_{n:n}$ represent their order statistics. For a fixed $p\in (0,1)$, define the binomial random variable $M_p={\rm card}\{i:U_i\leq p\}$, where ${\rm card}\{\cdot\}$ denotes the cardinality of a set.
Then, since the joint distribution of the order statistics $ X_{1:n},\dots,X_{n:n}$  is identical to that of $F^{-1}(U_{1:n}),\dots,F^{-1}(U_{n:n})$,
the remainder term $R_n$ in \eqref{eqn:residual} can be written as
\begin{equation}
\label{pr_1}
R_n\stackrel{d}{=} R_{n,1}+R_{n,2},
\end{equation}
where
\begin{align}
\label{pr_2}
\notag
R_{n,1}:=&\frac{np-[np]}{n(1-p)}\big(F^{-1}(U_{[np]+1:n})-F^{-1}(p)\big),\\
R_{n,2}:=&- \frac {{\rm {\rm {\rm sign}}}(M_p-[np])}{n(1-p)}\sum_{i=([np]\wedge M_p)+1}^{[np]\vee M_p} \big(F^{-1}(U_{i:n})-F^{-1}(p)\big).
\end{align}
Accordingly, we have
\begin{equation}
\label{pr_3}
-B_n=\mathbf{E}(R_n)=\mathbf{E}(R_{n,1})+\mathbf{E}(R_{n,2}).
\end{equation}

Further, let $\mathbb{U}_{p}\subset (0,1)$ denote a neighborhood of the point $p\in(0,1)$ in which the quantile function $F^{-1}$ is continuously differentiable (for the proof of Theorem~\ref{Thm_1}), or H\"older continuous of order $\gamma \in (0,1]$ (for the proof of Theorem~\ref{Thm_2}). In both cases, the following inequality holds \citep[e.g.,][pp.~453--454]{sw1986}:
\begin{equation}
\label{S_W_7}
\mathbf{P}\big( U_{[np]+1:n} \notin \mathbb{U}_p \big) \le e^{-c\, n},
\end{equation}
for some constant $c>0$.  The inequality above will be used repeatedly throughout the proofs of the two theorems.

\begin{proof}[\bf Proof of Theorem \ref{Thm_1}]  
Recall the relationship noted in \eqref{pr_3}.  In what follows,
we will prove that
\begin{equation}
\label{t1_01}
\mathbf{E}(R_{n,1})=o\big(n^{-1}\big)
\end{equation}
and that
\begin{equation}
\label{t1_02}
\mathbf{E}(R_{n,2})=\frac{p}{2nf(\xi_p)}+o\big(n^{-1}\big), 
\end{equation}
which together yield the desired asymptotic expansion in \eqref{B_n}. 

To do so, define the following shorthand notation:
$
\label{g(u)}
g(u):=\frac{\rd }{\rd u}\big(F^{-1}(u)\big)={1}/{f\big(F^{-1}(u)\big)}$ for $u\in\mathbb{U}_p.$
We fix $\epsilon_1>0$ and  $\delta>0$ such that $I_{\delta}:=(p-\delta,p+\delta)\subset \mathbb{U}_p$ and $|g(u)-g(v)|\leq \epsilon_1$ for all $u,v\in I_{\delta}$. 
To prove~\eqref{t1_01}, we write
\begin{equation}
\label{t1_1}
\mathbf{E}(R_{n,1})=\mathbf{E}_1+\mathbf{E}_2,
\end{equation}
where
\[
\mathbf{E}_1:=\mathbf{E}\big[R_{n,1}\mathds{1}_{ I_{\delta}}(U_{[np]+1:n})\big],\qquad \mathbf{E}_2:=\mathbf{E}\big[R_{n,1}\mathds{1}_{[0,1]\setminus I_{\delta}}(U_{[np]+1:n})\big].
\]
Hereafter, $\mathds{1}_D(\cdot)$ denotes the indicator function of the set $D$.
Let us first study $\mathbf{E}_1$. The following relationships hold: 
\begin{align}
\label{t1_2}
\notag
\mathbf{E}_1&=\frac{np-[np]}{n(1-p)}\mathbf{E}\Big[g(\zeta)\big(U_{[np]+1:n}-p\big)
\mathds{1}_{I_{\delta}}(U_{[np]+1:n})\Big]\\
&=\frac{np-[np]}{n(1-p)}\mathbf{E}\Big[\Big(g(p)\big(U_{[np]+1:n}-p\big)
+\left(g(\zeta)-g(p)\right)\big(U_{[np]+1:n}-p\big)\Big)
\mathds{1}_{I_{\delta}}(U_{[np]+1:n})\Big],
\end{align}
where $\zeta$ is a point lying between $p$ and $U_{[np]+1:n}$. 

Moreover, consider the following decomposition:  
\begin{equation}
\label{t1_3}
\mathbf{E}\Big[\Big(g(p)\big(U_{[np]+1:n}-p\big)
+\left(g(\zeta)-g(p)\right)\big(U_{[np]+1:n}-p\big)\Big)
\mathds{1}_{I_{\delta}}(U_{[np]+1:n})\Big]=e_1+e_2+e_3.
\end{equation}
For the first component in the above decomposition, we have
\begin{equation*}
\label{t1_4_1}
e_1:=g(p)\mathbf{E}\big(U_{[np]+1:n}-p\big)=g(p)\left(\frac{[np]+1
}{n+1}-p\right),
\end{equation*}
thus it holds that
\begin{equation}
\label{t1_4}
|e_1|\leq g(p)\,\frac{\max(p, 1-p)}{n+1}.
\end{equation}

For the second component which is given by
\begin{equation*}
\label{t1_5}
e_2:=- g(p)\mathbf{E}\Big[\big(U_{[np]+1:n}-p\big)(1-\mathds{1}_{I_{\delta}}(U_{[np]+1:n}))\Big],
\end{equation*}
we have
\begin{equation}
\label{t1_6}
|e_2|\leq g(p)\mathbf{E}\Big[1-\mathds{1}_{I_{\delta}}(U_{[np]+1:n})\Big]=g(p)\mathbf{P}\Big(U_{[np]+1:n}\notin I_{\delta}\Big)\leq g(p) \exp(-c\,n),
\end{equation}
where the rightmost inequality holds due to \eqref{S_W_7}.
Here and throughout, $c$ denotes a positive constant that does not depend on $n$ but may vary from line to line. 

Now, consider the third component in the decomposition \eqref{t1_3}, which is given by 
\begin{equation}
\label{t1_7}
e_3
:=\mathbf{E}\Big[\big(g(\zeta)-g(p)\big)\big(U_{[np]+1:n}-p\big)
\mathds{1}_{I_{\delta}}\big(U_{[np]+1:n}\big)\Big].
\end{equation}
By the choice of $\epsilon_1$ and $\delta$ noted at the beginning of this proof, and using well-known formulas for the variance of order statistics \citep[see, e.g.,][]{sw1986}, we can obtain
\begin{equation}
\label{t1_8}
|e_3|
\le \epsilon_1\,\mathbf{E}\big|U_{[np]+1:n}-p\big|
\le \epsilon_1\left(\frac{p(1-p)}{n}\big(1+O(n^{-1})\big)\right)^{1/2}.
\end{equation}
Combining \eqref{t1_2}–\eqref{t1_8} then yields
\begin{align}
\label{eqn:E1_result}
\mathbf{E}_1=o(n^{-1}).
\end{align}

Next, we turn to $\mathbf{E}_2$ in \eqref{t1_1}, which can be expressed as
\begin{equation*}
\label{t1_9}
\frac{np-[np]}{n(1-p)}\mathbf{E}\Big[\big(F^{-1}(U_{[np]+1:n})-F^{-1}(p)\big)\mathds{1}_{[0,1]\setminus I_{\delta}}(U_{[np]+1:n})\Big].
\end{equation*}
Let us consider the following decomposition: 
\begin{equation}
\label{t1_10}
\mathbf{E}\Big[\big(F^{-1}(U_{[np]+1:n})-F^{-1}(p)\big)\mathds{1}_{[0,1]\setminus I_{\delta}}(U_{[np]+1:n})\Big]=e_4+e_5,
\end{equation}
where
\begin{align*}
e_4&:=\mathbf{E}\big[F^{-1}(U_{[np]+1:n})\mathds{1}_{[0,1]\setminus I_{\delta}}(U_{[np]+1:n})\big],\\
e_5&:=-F^{-1}(p)\mathbf{E}\big[\mathds{1}_{[0,1]\setminus I_{\delta}}(U_{[np]+1:n})\big].
\end{align*}
For $e_4$,  it holds that
\begin{align}
\label{t1_11}
\notag
|e_4|&\leq \mathbf{E}\big[\big|F^{-1}(U_{[np]+1:n})\big|\mathds{1}_{[0,1]\setminus I_{\delta}}(U_{[np]+1:n})\big]\\
\notag
&\leq \big[\mathbf{E}\big|F^{-1}(U_{[np]+1:n})\big|^{1+\nu}\big]^{\frac{1}{1+\nu}}
\big[\mathbf{E}\big(\mathds{1}_{[0,1]\setminus I_{\delta}}(U_{[np]+1:n})\big)\big]^{1-\frac{1}{1+\nu}}\\
\notag
&\leq \underbrace{C(\rho)^{\frac{1}{1+\nu}}}_{C_1(\rho,\nu)}\left[\left(\mathbf{E}\big|X\big|^{\varepsilon} \frac{(n+1)^2}{([np]+1)(n-[np])}\right)^{\frac{1+\nu}{\varepsilon}}\right]^{\frac{1}{1+\nu}}\Big[ \mathbf{P}\Big(U_{[np]+1:n}\notin I_{\delta}\Big)\Big]^{1-\frac{1}{1+\nu}}\\
\notag
&\leq C_1(\rho,\nu)\left(\mathbf{E}\big|X\big|^{\varepsilon} \frac{(n+1)^2}{([np]+1)(n-[np])}\right)^{\frac{1}{\varepsilon}}\,\exp\left(-n\,\frac{c\nu}{1+\nu}\right)\\
&= C_1(\rho,\nu)\left(\frac{\mathbf{E}\big|X\big|^{\varepsilon}}{p(1-p)} \Big(1+O(n^{-1})\Big)\right)^{\frac{1}{\varepsilon}}\,\exp\left(-n\,\frac{c\nu}{1+\nu}\right).
\end{align}
In the above calculations, the second inequality follows from the H\"{o}lder inequality with an arbitrary $\nu>0$, while the third inequality applies the inequality in \eqref{L_es_mom1} with $\rho=(1+\nu)/\varepsilon$. 

For $c_5$, we get
\begin{equation}
\label{t1_12}
|e_5|=\big|F^{-1}(p)\big|\,\mathbf{P}\Big(U_{[np]+1:n}\notin I_{\delta}\Big)\leq \big|F^{-1}(p)\big|\,\exp(-c\, n).
\end{equation}
Combining \eqref{t1_10}--\eqref{t1_12} then implies
\begin{align}
\label{eqn:E2_result}
\mathbf{E}_2=o(n^{-1}).
\end{align}
Collectively, the relationship in \eqref{t1_01} holds because of \eqref{eqn:E1_result} and \eqref{eqn:E2_result}.

So far, we have established \eqref{t1_01}. Now, we proceed to proving \eqref{t1_02}, which is more demanding. 
Recall that $\epsilon_1$ and $\delta$ are defined and fixed at the beginning of the proof of this theorem. We note that if $U_{[np]+1:n}\in I_{\delta}$, then  $U_{i:n}\in I_{\delta}$ for all $i=([np]\wedge M_p)+1,\dots,[np]\vee M_p$. Indeed, let $U_{[np]+1:n}\in I_{\delta}$. If $M_p>[np]$, then  $U_{[np]+1:n}\leq U_{M_p:n}\leq p$ and $(U_{[np]+1:n},U_{M_p:n})\subset I_{\delta}$. If $M_p<[np]$, then $p<U_{M_p+1:n}\leq U_{[np]+1:n}$, and hence, $(U_{M_p+1:n},U_{[np]:n})\subset I_{\delta}$. If $M_p=[np]$, then $R_{n,2}$ is zero. In passing, we note that, whenever $M_p\neq [np]$, the following inequality holds:
\begin{equation}
\label{t1_for_max}
\max_{([np]\wedge M_p)+1 \le i \le ([np]\vee M_p)}
\left| F^{-1}(U_{i:n}) - F^{-1}(p) \right|
\le
\left| F^{-1}(U_{[np]:n}) - F^{-1}(p) \right|.
\end{equation}
This inequality will be useful later.

We evaluate $\mathbf{E}(R_{n,2})$ in \eqref{t1_02} by decomposing it as
\begin{equation}
\label{t1_13}
\mathbf{E}(R_{n,2})=\mathbf{E}_3+\mathbf{E}_4,
\end{equation}
where
\[
\mathbf{E}_3:=\mathbf{E}\big[R_{n,2}\mathds{1}_{ I_{\delta}}(U_{[np]+1:n})\big],\qquad \mathbf{E}_4:=\mathbf{E}\big[R_{n,2}\mathds{1}_{[0,1]\setminus I_{\delta}}(U_{[np]+1:n})\big].
\]
We first treat $\mathbf{E}_3$. Consider $R_{n,2}$ on the event $\big\{U_{[np]+1:n}\in I_{\delta}\big\}$,  on which
$\mathds{1}_{ I_{\delta}}(U_{[np]+1:n})=1$. On this event, we have
\begin{align}
\label{t1_14}
\notag
R_{n,2}=-& \frac {{\rm {\rm {\rm sign}}}(M_p-[np])}{n(1-p)}\sum_{i=([np]\wedge M_p)+1}^{[np]\vee M_p} g(p) \big(U_{i:n}-p\big)\\
&\qquad -\frac {{\rm {\rm {\rm sign}}}(M_p-[np])}{n(1-p)}\sum_{i=([np]\wedge M_p)+1}^{[np]\vee M_p}\big( g(\zeta_i)-g(p) \big) \big(U_{i:n}-p\big)=:R'_{n,2}+R''_{n,2},
\end{align}
where each $\zeta_i$ lies between $p$ and $U_{i:n}$, $R'_{n,2}$ and $R''_{n,2}$ denote the first and second summation terms in \eqref{t1_14}, respectively. Accordingly, $\mathbf{E}_3$ can be written as 
\begin{align}
\label{t1_15}
\notag
\mathbf{E}_3&=\mathbf{E}\big[\big(R'_{n,2}+R''_{n,2}\big)\mathds{1}_{ I_{\delta}}(U_{[np]+1:n})\big]\\
&=
\mathbf{E}\big(R'_{n,2}\big) -
\mathbf{E}\big[R'_{n,2}\big(1-\mathds{1}_{ I_{\delta}}(U_{[np]+1:n})\big)\big]+\mathbf{E}\big[R''_{n,2}\mathds{1}_{ I_{\delta}}(U_{[np]+1:n})\big]\notag\\
&=:\mathbf{E}_3'+e_6+e_7.
\end{align}
We claim that $|e_6+e_7|=o\big(n^{-1}\big)$.  To see this, for $e_6$, we have
\begin{align*}
\notag
|e_6|&\leq \frac{g(p)}{n (1-p)}\mathbf{E}\Big(|M_p-[np]|\,
\big(1-\mathds{1}_{ I_{\delta}}(U_{[np]+1:n})\big)
\Big)\\
&\leq \frac{g(p)}{n(1-p)}\Big(\mathbf{E}(|M_p-[np]|)^2\mathbf{P}(U_{[np]+1:n}\notin I_{\delta})\Big)^{1/2}
\leq \frac{g(p)p^{1/2}}{(n(1-p))^{1/2}} \exp(-c\, n).
\end{align*}
For $e_7$, we get
\begin{align*}
\label{t1_17}
\notag
|e_7|&\leq \frac{\sup_{u\in I_{\delta}}|g(u)- g(p)|}{n (1-p)}\mathbf{E}\Big(\big|M_p-[np]\big|\,
\big|U_{[np]:n}-p\big|\Big)\\
\notag
&\leq \frac{\epsilon_1}{n (1-p)}\Big(\mathbf{E}\big(M_p-[np]\big)^2\, \mathbf{E}\big(U_{[np]:n}-p\big)^2 \Big)^{1/2}\\
&=\frac{\epsilon_1\big(np(1-p)\big)^{1/2}}{n (1-p)}\left(\frac{p(1-p)}{n}\big(1+O(n^{-1})\big)\right)^{1/2}=\frac{\epsilon_1 p}{n}+O\big(n^{-2}\big),
\end{align*}
where we used the inequality in \eqref{t1_for_max} and the Cauchy--Schwarz inequality to establish the desired relationships. Since $\epsilon_1>0$ is arbitrary, it follows that $e_7=o(n^{-1})$.

We now proceed to evaluating the first term in \eqref{t1_15}. We note that the following relationships hold:
\begin{equation}
\label{t1_18}
\mathbf{E}_3'=\mathbf{E}\big(R'_{n,2}\big)
=-\frac{g(p)}{n(1-p)}\mathbf{E}\big(\Sigma_n\big),
\end{equation}
where
\begin{align*}
\label{l3_1}
\Sigma_n&:={\rm {\rm {\rm sign}}}(M_p-[np])\sum_{i=([np]\wedge M_p)+1}^{[np]\vee M_p} \big(U_{i:n}-p\big)=\Sigma'_n+\Sigma''_n,\\
\notag
\Sigma'_n&:=\mathds{1}_{\{M_p>[np]\}}\sum_{i=[np]+1}^{M_p} \big(U_{i:n}-p\big)=T_{1,n}-\mathds{1}_{\{M_p>[np]\}} (M_p-[np])\big(p-U_{M_p:n}\big),\\
\notag
\Sigma''_n&:=-
\mathds{1}_{\{M_p<[np]\}}\sum_{i=M_p+1}^{[np]} \big(U_{i:n}-p\big)
=T_{2,n}-\mathds{1}_{\{M_p<[np]\}} ([np]-M_p)\big(U_{M_p+1:n}-p\big),
\end{align*}
with
\[
T_{1,n}:=-\mathds{1}_{\{M_p>[np]\}}\sum_{i=[np]+1}^{M_p} \big(U_{M_p:n}-U_{i:n}\big),\quad
T_{2,n}:=-\mathds{1}_{\{M_p<[np]\}}\sum_{i=M_p+1}^{[np]} \big(U_{i:n}-U_{M_p+1:n}\big).
\]
Moreover, note that if  $M_p=[np]+1$, then  $T_{1,n}=0$. Therefore, we can slightly simplify $T_{1,n}$ into
\begin{equation*}
\label{l3_4}
T_{1,n}=-\mathds{1}_{\{M_p>[np]+1\}}\sum_{i=[np]+1}^{M_p-1} \big(U_{M_p:n}-U_{i:n}\big)=-\mathds{1}_{\{M_p>[np]+1\}}\sum_{i=[np]+1}^{M_p-1}\sum_{j=i+1}^{M_p}(U_{j:n}-U_{j-1:n}).
\end{equation*}
Similarly, if  $M_p=[np]-1$, then $T_{2,n}=0$.  Thus, $T_{2,n}$ can be simplified into 
\begin{equation*}
\label{l3_5}
T_{2,n}=-\mathds{1}_{\{M_p<[np]-1\}}\sum_{i=M_p+2}^{[np]} \big(U_{i:n}-U_{M_p+1:n}\big)=-\mathds{1}_{\{M_p<[np]-1\}}\sum_{i=M_p+2}^{[np]}\sum_{j=M_p+2}^{i}(U_{j:n}-U_{j-1:n}).
\end{equation*}
Let $s_{j,n} := U_{j:n} - U_{j-1:n}$, $j = 1,\ldots,n+1$, which correspond to the spacings of the uniform order statistics on $[0,1]$, with $U_{0:n}=0$ and $U_{n+1:n}=1$.  Then we can rewrite $T_{1,n}$ and $T_{2,n}$ as 
\begin{align*}
\notag
T_{1,n}=-\mathds{1}_{\{M_p>[np]+1\}}\sum_{i=[np]+1}^{M_p-1}\sum_{j=i+1}^{M_p}s_{j,n}=
&-\mathds{1}_{\{M_p>[np]+1\}}\sum_{i=[np]+1}^{M_p-1}\left(\sum_{j=i+1}^{M_p}s_{j,n}-\frac{M_p-i}{n+1}\right)\\
\notag
&-\mathds{1}_{\{M_p>[np]+1\}}\frac{(M_p-[np]-1)(M_p-[np])}{2(n+1)},\\
\notag
T_{2,n}=-\mathds{1}_{\{M_p<[np]-1\}}\sum_{i=M_p+2}^{[np]}\sum_{j=M_p+2}^{i}s_{j,n}=
&-\mathds{1}_{\{M_p<[np]-1\}}\sum_{i=M_p+2}^{[np]}\left(\sum_{j=M_p+2}^{i}s_{j,n}-\frac{i-M_p-1}{n+1}\right)\\
&-\mathds{1}_{\{M_p<[np]-1\}}\frac{([np]-M_p-1)([np]-M_p)}{2(n+1)}.
\end{align*}

In addition, we observe that
\begin{align}
\label{indicators-2}
-\mathds{1}_{\{M_p>[np]+1\}}\frac{(M_p-[np]-1)(M_p-[np])}{2(n+1)}-\mathds{1}_{\{M_p<[np]-1\}}\frac{([np]-M_p-1)([np]-M_p)}{2(n+1)}
\end{align}
can be written equivalently as
\begin{equation}
\label{vervaat}
- \frac{(|M_p-[np]|-1)|M_p-[np]|}{2(n+1)}=:V_{p,n}.
\end{equation}
Indeed, the expression in \eqref{vervaat} equals zero when
\( M_p=[np] \), \( M_p=[np]+1 \), or \( M_p=[np]-1 \), which matches exactly the cases excluded by the indicator functions involved in \eqref{indicators-2}.

Collectively, we can conclude
\begin{equation}
\label{l3_7}
\Sigma_n=V_{p,n}+r_{1,n}+r_{2,n}+r_{3,n},
\end{equation}
where
\begin{align*}
r_{1,n}&:=-\mathds{1}_{\{M_p>[np]\}} (M_p-[np])\big(p-U_{M_p:n}\big)-\mathds{1}_{\{M_p<[np]\}} ([np]-M_p)\big(U_{M_p+1:n}-p\big),\\
r_{2,n}&:=-\mathds{1}_{\{M_p>[np]+1\}}\sum_{i=[np]+1}^{M_p-1}\left(\sum_{j=i+1}^{M_p}s_{j,n}-\frac{M_p-i}{n+1}\right),\\
r_{3,n}&:=-\mathds{1}_{\{M_p<[np]-1\}}\sum_{i=M_p+2}^{[np]}\left(\sum_{j=M_p+2}^{i}s_{j,n}-\frac{i-M_p-1}{n+1}\right).
\end{align*}

For the first term on the right-hand side of \eqref{l3_7}, we have
\begin{align}
\label{l3_8}
\notag
\mathbf{E}\big(V_{p,n}\big)&=-\frac{1}{2(n+1)}\mathbf{E}\Big((M_p-[np])^2-|M_p-[np]|\Big)\\
\notag
&=-\frac{1}{2(n+1)}\mathbf{E}\Big((M_p-np)^2+(np-[np])(M_p-[np]+M_p-np) -|M_p-[np]|\Big)\\
&=-\frac{np(1-p)}{2(n+1)}+\delta_n,
\end{align}
where
\begin{align*}
\label{l3_9}
\delta_n:=-\frac{1}{2(n+1)} \mathbf{E}\Big((np-[np])(M_p-[np]+M_p-np) -|M_p-[np]|\Big).
\end{align*}
Since $np-[np]<1$ and $|M_p-[np]|<|M_p-np|+1$, it follows that
\begin{equation}
\label{l3_10}
|\delta_n|<\frac{3\mathbf{E}|M_p-np|+2}{2(n+1)} < \frac{3\big(\mathbf{E}(M_p-np)^2\big)^{1/2}+2}{2n}=\frac{3\big(p(1-p)\big)^{1/2}}{2n^{1/2}}<\frac{3}{4n^{1/2}}<\frac{1}{n^{1/2}}.
\end{equation}

Next we develop bounds for $r_{k,n}$,  $k=1,2,3$, appearing on the right-hand side of~\eqref{l3_7}. For $r_{1,n}$, we have
\begin{align}
\label{l3_11}
\notag
|r_{1,n}|&\leq |M_p-[np]\,|(p-U_{M_p:n})+|M_p-[np]\,|(U_{M_p+1:n}-p)=|M_p-[np]\,|(U_{M_p+1:n}-U_{M_p:n})\\
&\leq \big(|M_p-np|+1\big)(U_{M_p+1:n}-U_{M_p:n})\stackrel{d}{=}\big(|M_p-np|+1\big)U_{1:n}.
\end{align}
where, in the last step, we used the fact that
$(U_{M_p+1:n}-U_{M_p:n}) = s_{M_p+1,n} \stackrel{d}{=} s_{1,n} = U_{1:n}$; 
that is, the uniform spacings are identically distributed and exchangeable random variables
\cite[see,][p.~721]{sw1986}.
The relationships derived in~\eqref{l3_11} imply that
\begin{align}
\label{l3_12}
\notag
|\mathbf{E}(r_{1,n})|&\leq \mathbf{E}\big(|M_p-np|\,U_{1:n}\big)+\mathbf{E}U_{1:n}\leq \Big(np(1-p)\mathbf{E}\big(U_{1:n}\big)^2\Big)^{1/2}+n^{-1}\\
&=\left(\frac{2np(1-p)}{(n+2)(n+1)}\right)^{1/2}+n^{-1}<(2n)^{-1/2}+n^{-1}<2n^{-1/2}.
\end{align}

We now turn to the terms $r_{2,n}$ and $r_{3,n}$ involved in \eqref{l3_7}. Let us set $k=i-[np]$, then $r_{2,n}$ can be rewritten as
\begin{align}
\label{l3_r2}
\notag
r_{2,n}&=-\mathds{1}_{\{M_p>[np]+1\}}\sum_{k=1}^{M_p-[np]-1}\left(\sum_{j=[np]+k+1}^{M_p}s_{j,n}-\frac{M_p-[np]-k}{n+1}\right)\\
\notag
&\stackrel{d}{=}-\mathds{1}_{\{M_p>[np]+1\}}\sum_{k=1}^{M_p-[np]-1}\left(\sum_{j=1}^{M_p-[np]-k}s_{j,n}-\frac{M_p-[np]-k}{n+1}\right)\\
\notag
&=-\mathds{1}_{\{M_p>[np]+1\}}\sum_{k=1}^{M_p-[np]-1}\left(U_{M_p-[np]-k:n}-\frac{M_p-[np]-k}{n+1}\right)\\
&= -\mathds{1}_{\{M_p>[np]+1\}}\sum_{i=1}^{M_p-[np]-1}\left(U_{i:n}-\frac{i}{n+1}\right),
\end{align}
where in the second line, we used again the fact that the uniform spacings are identically distributed and exchangeable random variables.
 
Similarly, setting $k=i-M_p-1$, we can obtain the following relationships for $r_{3,n}$:
\begin{align}
\label{l3_r3}
\notag
r_{3,n}&=-\mathds{1}_{\{M_p<[np]-1\}}\sum_{k=1}^{[np]-M_p-1}\left(\sum_{j=M_p+2}^{k+M_p+1}s_{j,n}-\frac{k}{n+1}\right)\\
\notag
&\stackrel{d}{=}-\mathds{1}_{\{M_p<[np]-1\}}\sum_{k=1}^{[np]-M_p-1}\left(\sum_{j=1}^{k}s_{j,n}-\frac{k}{n+1}\right)\\
&=-\mathds{1}_{\{M_p<[np]-1\}}\sum_{k=1}^{[np]-M_p-1}\left(U_{k:n}-\frac{k}{n+1}\right).
\end{align} 
Then combining \eqref{l3_r2} and \eqref{l3_r3} yields
\begin{align}
\label{l3_14}
\notag
&|\mathbf{E}(r_{2,n}+r_{3,n})|\leq \mathbf{E}|r_{2,n}|+\mathbf{E}|r_{3,n}|\\
\notag
\leq & \mathbf{E}\left( \mathds{1}_{\{M_p>[np]+1\}}\sum_{i=1}^{M_p-[np]-1}\left|U_{i:n}-\frac{i}{n+1}\right|  \right) +
\mathbf{E}\left( \mathds{1}_{\{M_p<[np]-1\}}\sum_{i=1}^{[np]-M_p-1}\left|U_{i:n}-\frac{i}{n+1}\right|  \right)\\
\notag
=& \mathbf{E}\left( \mathds{1}_{\{M_p>[np]+1\}}\sum_{i=1}^{|M_p-[np]-1\,|}\left|U_{i:n}-\frac{i}{n+1}\right|   +
\mathds{1}_{\{M_p<[np]-1\}}\sum_{i=1}^{|\,[np]-M_p-1|}\left|U_{i:n}-\frac{i}{n+1}\right|  \right)\\
\leq &\mathbf{E}\left(\sum_{i=1}^{|M_p-[np]| +1}\left|U_{i:n}-\frac{i}{n+1}\right| \right). 
\end{align}

Let $A>1$ be arbitrary. By Hoeffding’s inequality \cite[][Theorem~1]{Hoeffding1963}, we know
\[
\mathbf{P}\left(|M_p-np|\geq \sqrt{A\,n\log n}\right)\leq 2\exp\left(-2A\log n\right)=2n^{-2A}<2n^{-2}
\]
for all sufficiently large $n$ such that $0<\sqrt{A\log n / n}<p$.
Therefore, the right-hand side of \eqref{l3_14} can be bounded by
\begin{align}
\label{t1_26}
\notag
&\sum_{i=1}^{\lceil\sqrt{A\,n\log n} \rceil}\mathbf{E}\left|U_{i:n}-\frac{i}{n+1}\right|+      \mathbf{E}\left(\mathds{1}_{\{|M_p-np|\geq \sqrt{A\,n\log n}\}}\sum_{i=1}^{|M_p-[np]\,|+1}\left|U_{i:n}-\frac{i}{n+1}\right|\right)\\
\notag
\leq &\sum_{i=1}^{\lceil\sqrt{A\,n\log n} \rceil}
\left(\frac{i(n-i+1)}{(n+1)^3} \right)^{1/2} +n\,\mathbf{P}\left(|M_p-np|\geq \sqrt{A\,n\log n}\right)\\
\notag
\leq &\sum_{i=1}^{\lceil\sqrt{A\,n\log n} \rceil}
\left(\frac{i}{(n+1)^2} \right)^{1/2} +2n^{-1}\leq \frac{\lceil\sqrt{A\,n\log n} \rceil}{n+1}(\lceil\sqrt{A\,n\log n} \rceil)^{1/2}+2n^{-1}\\
=&\frac{(A\,n\log n)^{3/4}}{n}+O(n^{-1})=O\left(\frac{(\log n)^{3/4}}{ n^{1/4}}\right),
\end{align}
where $\lceil\cdot\rceil$ denotes the ceiling function.

Collecting the relationships established in \eqref{t1_18}, \eqref{l3_7}--\eqref{l3_8} and the estimates
\eqref{l3_10}, \eqref{l3_12}, and \eqref{t1_26}, we readily obtain
\begin{align*}
\label{t1_27}
\notag
\mathbf{E}_3' &=-\frac{g(p)}{n(1-p)}\mathbf{E}\left[V_{p,n} +r_{1,n}+r_{2,n}+r_{3,n}\right] 
=\frac{g(p)}{n(1-p)}\left[\frac{np(1-p)}{2(n+1)}
+O\left(\frac{(\log n)^{3/4}}{ n^{1/4}}\right)\right]\\
&=\frac{g(p)p}{2(n+1)}
+O\left(\frac{(\log n)^{3/4}}{n^{5/4}}\right)=\frac{p}{2nf(\xi_p)}+o\big(n^{-1}\big).
\end{align*}

To complete the proof, it remains to show that the term $\mathbf{E}_4$ in \eqref{t1_13} satisfies
\begin{equation}
\label{t1_29}
\mathbf{E}_4=o\big(n^{-1}\big).
\end{equation} 
To this end, the following bound is established: 
\begin{align}
\label{t1_30}
\notag
|\mathbf{E}_4| &\leq \mathbf{E}\big[|R_{n,2}|\mathds{1}_{[0,1]\setminus I_{\delta}}(U_{[np]+1:n})\big]\\
\notag
&= \frac {1}{n(1-p)}\mathbf{E}\left[\left|
\sum_{i=([np]\wedge M_p)+1}^{[np]\vee M_p} \big(F^{-1}(U_{i:n})-F^{-1}(p)\big)\right|\mathds{1}_{[0,1]\setminus I_{\delta}}(U_{[np]+1:n})\right]\\
\notag
&\leq  \frac {1}{1-p}\mathbf{E}\left[\left|F^{-1}(U_{[np]:n})-F^{-1}(p)\right|\mathds{1}_{[0,1]\setminus I_{\delta}}(U_{[np]+1:n})\right]\\
&\leq  \frac {1}{1-p}\mathbf{E}\left[\left|F^{-1}(U_{[np]:n})\right|\mathds{1}_{[0,1]\setminus I_{\delta}}(U_{[np]+1:n})\right]+\frac {\big|F^{-1}(p)\big|}{1-p}\mathbf{P}\Big(U_{[np]+1:n}\notin I_{\delta}\Big),
\end{align}
where, in particular, the inequality in \eqref{t1_for_max} was used in the third step.
The first term on the right-hand side of \eqref{t1_30} has already been bounded in \eqref{t1_11}. Consequently, for an arbitrary $\nu>0$, the right-hand side of \eqref{t1_30} is bounded above by
\begin{align*}
\frac { C_1(\rho,\nu)}{1-p}\left(\frac{\mathbf{E}\big|X\big|^{\varepsilon}}{p(1-p)} \Big(1+O(n^{-1})\Big)\right)^{\frac{1}{\varepsilon}}\,\exp\left(-n\,\frac{c\nu}{1+\nu}\right)
+\frac{\big|F^{-1}(p)\big|}{1-p}\exp(-c\,n),
\end{align*}
where $C_1(\rho,\nu)>0$ is a constant depending only on $\nu$ and $\rho=(1+\nu)/\varepsilon$, as defined in~\eqref{t1_11}. 
Since both terms in the above expression decay exponentially fast in $n$, this bound implies \eqref{t1_29}.

The proof of Theorem~\ref{Thm_1} is finally completed.
\end{proof}

\bigskip

\begin{proof}[{\bf Proof of Theorem~\ref{Thm_2}}] Recall that, from \eqref{pr_1}--\eqref{pr_3}, we have 
\begin{equation}
\label{t2_1}
-B_n=\mathbf{E}(R_n)=\mathbf{E}(R_{n,1})+\mathbf{E}(R_{n,2}).
\end{equation} 
We shall first show that $\mathbf{E}(R_{n,1})$ is of negligible order relative to the right-hand side of \eqref{B_n_estim}. 
To this end, we write
\begin{equation*}
\label{t2_2}
\mathbf{E}(R_{n,1})=\mathbf{E}_1+\mathbf{E}_2,
\end{equation*}
where
\begin{equation*}
\label{t2_3}
\mathbf{E}_1:=\mathbf{E}\big(R_{n,1}\mathds{1}_{\mathbb{U}_p}(U_{[np]+1:n})\big),\qquad \mathbf{E}_2:=\mathbf{E}\big(R_{n,1}\mathds{1}_{[0,1]\setminus\mathbb{U}_p}(U_{[np]+1:n})\big). 
\end{equation*}
Here, $\mathbb{U}_p$ denotes the neighborhood of $p$ on which $F^{-1}$ is locally H\"older continuous of order $\gamma\in(0,1]$ with constant $c_\gamma$, as specified in \eqref{Lipschitz}.

For $\mathbf{E}_1$, we get
\begin{align*}
\notag
|\mathbf{E}_1|&\leq \frac{c_{\gamma}}{n(1-p)}\mathbf{E}|U_{[np]+1:n}-p|^{\gamma}\leq
\frac{c_{\gamma}}{n(1-p)}\Big(\mathbf{E}(U_{[np]+1:n}-p)^2\Big)^{\gamma/2}\\
&=\frac{c_{\gamma}}{n(1-p)}\Big(\frac{p(1-p)}{n}(1+o(1))\Big)^{\gamma/2}=
\frac{c_{\gamma}p^{\gamma/2}}{n^{1+\gamma/2}(1-p)^{1-\gamma/2}}\big(1+o(1)\big).
\end{align*}
Therefore, $|\mathbf{E}_1|=O(n^{-1-\gamma/2})=o(n^{-1})$, which is negligible relative to the bound in \eqref{B_n_estim}.


For $\mathbf{E}_2$, we have
\begin{align*}
\label{t2_5}
\notag
|\mathbf{E}_2|&=\frac{np-[np]}{n(1-p)}\left|
\mathbf{E}\Big(\big(F^{-1}(U_{[np]+1:n})-F^{-1}(p)\big)\mathds{1}_{[0,1]\setminus\mathbb{U}_p}(U_{[np]+1:n})\Big)\right|\\
&\leq \frac{1}{n(1-p)}\left(\mathbf{E}\Big|F^{-1}(U_{[np]+1:n})\mathds{1}_{[0,1]\setminus\mathbb{U}_p}(U_{[np]+1:n})\Big|
+|F^{-1}(p)|\mathbf{P}\big(U_{[np]+1:n} \notin \mathbb{U}_p\big)\right).
\end{align*}
Arguing as in the estimation of $e_4$ in \eqref{t1_11}, it follows that for any $\nu>0$,
\begin{equation*}
\label{t2_6}
|\mathbf{E}_2|
\le\frac{1}{n(1-p)}\left(
C_1(\rho,\nu)\left(\frac{\mathbf{E}\big|X\big|^{\varepsilon}}{p(1-p)} \Big(1+O(n^{-1})\Big)\right)^{\frac{1}{\varepsilon}}\,\exp\left(-n\,\frac{c\nu}{1+\nu}\right)
+|F^{-1}(p)|e^{-c\,n}\right),
\end{equation*}
where $\rho=(1+\nu)/\varepsilon$.  Since both terms in the above expression decay exponentially fast in $n$, we can conclude that
$\mathbf{E}_2=o(n^{-1})$, and hence $\mathbf{E}_2$ is negligible relative to the bound in \eqref{B_n_estim}.

It remains to treat $\mathbf{E}(R_{n,2})$ in \eqref{t2_1}. We decompose it as
\begin{equation*}
\label{t2_7}
\mathbf{E}(R_{n,2})=\mathbf{E}_3+\mathbf{E}_4,
\end{equation*}
where
\begin{equation*}
\label{t2_8}
\mathbf{E}_3:=\mathbf{E}\big(R_{n,2}\mathds{1}_{\mathbb{U}_p}(U_{[np]+1:n})\big),\qquad \mathbf{E}_4:=\mathbf{E}\big(R_{n,2}\mathds{1}_{[0,1]\setminus\mathbb{U}_p}(U_{[np]+1:n})\big).
\end{equation*}
We first consider the term $\mathbf{E}_3$. 
Note that all terms in the sum defining $R_{n,2}$ in \eqref{pr_2} have the same sign, either all positive or all negative, so the absolute value of the sum equals the sum of the absolute values.
By the local H\"older continuity of $F^{-1}$, we can obtain
\begin{align}
\label{t2_9}
|\mathbf{E}_3|\leq \frac{c_{\gamma}}{n(1-p)}\mathbf{E}\left(
\sum_{i=([np]\wedge M_p)+1}^{[np]\vee M_p} \big|U_{i:n}-p\big|^{\gamma}\right).
\end{align}
Moreover, by \eqref{t1_for_max}, the largest term in the sum in \eqref{t2_9} is bounded above by $|U_{[np]:n}-p|^{\gamma}$. Consequently, the right-hand side of \eqref{t2_9} is bounded by
\begin{align}
\label{t2_10}
\notag
&\frac{c_{\gamma}}{n(1-p)}\mathbf{E}\Big(|M_p-[np]| |U_{[np]:n}-p|^{\gamma}\Big)\\
\notag
\leq &\frac{c_{\gamma}}{n(1-p)}\mathbf{E}\Big(\big(|M_p-np|+1\big) |U_{[np]:n}-p|^{\gamma}\Big)\\
\notag
\leq &\frac{c_{\gamma}}{n(1-p)}\Big[\mathbf{E}\big(|M_p-np|+1\big)^2 \,\mathbf{E} |U_{[np]:n}-p|^{2\gamma}\Big]^{1/2}\\
\notag
\leq
&\frac{c_{\gamma}}{n(1-p)}\Big[np(1-p)+2\mathbf{E}|M_p-np|+1\Big]^{1/2}\Big[\big(\mathbf{E}(U_{[np]:n}-p)^2\big)^{2\gamma/2}\Big]^{1/2} \\
\notag
\leq &\frac{c_{\gamma}}{n(1-p)}\Big[np(1-p)+2(np(1-p))^{1/2}+1\Big]^{1/2}\left[\frac{p(1-p)}{n}\big(1+O(n^{-1})\big)\right]^{\gamma/2} \\
\notag
=&\frac{c_{\gamma}}{n(1-p)}\Big[np(1-p)\big(1+O(n^{-1/2})\big)\Big]^{1/2} \left[\frac{p(1-p)}{n}\big(1+O(n^{-1})\big)\right]^{\gamma/2} \\
\notag
=&\frac{c_{\gamma}}{n(1-p)}\Big[np(1-p)\Big]^{1/2}\big(1+O(n^{-1/2})\big)\left[\frac{p(1-p)}{n}\right]^{\gamma/2}\big(1+O(n^{-1})\big)\\
= &\frac{c_{\gamma}p^{\frac{1+\gamma}{2}}}{n^{\frac{1+\gamma}{2}}(1-p)^{\frac{1-\gamma}{2}}}\big(1+O(n^{-1/2})\big).
\end{align}
The above result leads exactly to the expression of the bound \eqref{B_n_estim} in Theorem~\ref{Thm_2}.

To complete the proof, it remains to show that the term $\mathbf{E}_4$ is of negligible order relative to the bound on the right-hand side of \eqref{t2_10}. 
Proceeding as before, but now taking into account the indicator
$\mathds{1}_{[0,1]\setminus\mathbb{U}_p}(U_{[np]+1:n})$, we obtain
\begin{align*}
\label{t2_11}
\notag
|\mathbf{E}_4|&\leq \frac{1}{n(1-p)}\mathbf{E}\Big(|M_p-[np]| \big|F^{-1}(U_{[np]:n})-F^{-1}(p)\big|\mathds{1}_{[0,1]\setminus\mathbb{U}_p}(U_{[np]+1:n})\Big)\\
\notag
&\leq \frac{1}{1-p}\mathbf{E}\Big(\big|F^{-1}(U_{[np]:n})-F^{-1}(p)\big|  \mathds{1}_{[0,1]\setminus\mathbb{U}_p}(U_{[np]+1:n})\Big)\\
\notag
&\leq \frac{1}{1-p}\mathbf{E}\Big(\big|F^{-1}(U_{[np]:n})\big|  \mathds{1}_{[0,1]\setminus\mathbb{U}_p}(U_{[np]+1:n})\Big)+\frac{\big|F^{-1}(p)\big|}{1-p}\mathbf{P}\big(U_{[np]+1:n} \notin \mathbb{U}_p\big)\\
&\leq \frac{C_1(\rho,\nu)}{1-p}\left(\frac{\mathbf{E}\big|X\big|^{\varepsilon}}{p(1-p)} \Big(1+O(n^{-1})\Big)\right)^{\frac{1}{\varepsilon}}\exp\left(-n\,\frac{c\nu}{1+\nu}\right) +\frac{\big|F^{-1}(p)\big|}{1-p} \exp(-c\,n),
\end{align*}
for each $\nu>0$ and $\rho=(1+\nu)/\varepsilon$; see a similar derivation in \eqref{t1_11}. This implies that $\mathbf{E}_4$ is of exponentially small order and therefore negligible relative to the bound in \eqref{t2_10}. 
This completes the proof of Theorem~\ref{Thm_2}.
\end{proof}

\end{appendices}

\end{document}